\documentclass[final,3p,times]{elsarticle}

\usepackage{amssymb}
\usepackage{amsmath}
 \usepackage{amsthm}
  \usepackage{color}

\usepackage{subfigure}
\usepackage{hyperref}

\newtheorem{thm}{Theorem}[section]

\newtheorem{rem}[thm]{Remark}

\newtheorem{asm}[thm]{Assumption}

\def\ba{{\mathbf{a}}}
\def\cN{{\cal N}}  

\def\bE{{\bf E}}

\def\bW{{\bf W}}

\newcommand{\C}{\rm I\kern-.5emC}
\newcommand{\R}{\rm I\kern-.19emR}

\def\3bar{{|\hspace{-.02in}|\hspace{-.02in}|}}

\newcommand{\jump}[1]{[\![ #1]\!]}

\journal{Elsevier}

\begin{document}
\newcommand{\BX}{{\bf X}}
\newcommand{\cv}{{\cal V}}
\newcommand{\cW}{{\cal W}}
\newcommand{\co}{{\cal O}}

\renewcommand{\theequation}{\thesection.\arabic{equation}}
\def\@eqnnum{{\reset@font\rm (\theequation)}}

\def\Xint#1{\mathchoice
{\XXint\displaystyle\textstyle{#1}}%
{\XXint\textstyle\scriptstyle{#1}}%
{\XXint\scriptstyle\scriptscriptstyle{#1}}%
{\XXint\scriptscriptstyle\scriptscriptstyle{#1}}%
\!\int}
\def\XXint#1#2#3{{\setbox0=\hbox{$#1{#2#3}{\int}$}
\vcenter{\hbox{$#2#3$}}\kern-.5\wd0}}
\def\ddashint{\Xint=}
\def\dashint{\Xint-}

\def\a{\alpha}
\def\b{\beta}
\def\d{\delta}\def\D{\Delta}
\def\e{\epsilon}
\def\g{\gamma}\def\G{\Gamma}
\def\k{\kappa}
\def\lam{\lambda}\def\Lam{\Lambda}
\renewcommand\o{\omega}\renewcommand\O{\Omega}
\def\s{\sigma}\def\S{\Sigma}
\renewcommand\t{\theta}\def\vt{\vartheta}
\newcommand{\vphi}{\varphi}
\def\z{\zeta}

\newcommand{\tsigma}{\tilde{\s}}
\newcommand{\tbsigma}{\tilde{\bsigma}}
\def\te{\tilde{\e}}
\def\tu{\tilde{u}}

\newcommand{\bchi}{\mbox{\boldmath$\chi$}}
\newcommand{\bdelta}{\mbox{\boldmath$\delta$}}
\newcommand{\bepsilon}{\mbox{\boldmath$\epsilon$}}
\newcommand{\bfeta}{\mbox{\boldmath$\eta$}}
\newcommand{\bgamma}{\mbox{\boldmath$\gamma$}}
\newcommand{\bomega}{\mbox{\boldmath$\omega$}}
\newcommand{\bvphi}{\mbox{\boldmath$\varphi$}}
\newcommand{\bphi}{\mbox{\boldmath$\phi$}}
\newcommand{\bPhi}{\mbox{\boldmath$\Phi$}}
\newcommand{\bpsi}{\mbox{\boldmath$\psi$}}
\newcommand{\bPsi}{\mbox{\boldmath$\Psi$}}
\newcommand{\bsigma}{\mbox{\boldmath$\sigma$}}
\newcommand{\btau}{\mbox{\boldmath$\tau$}}
\newcommand{\bxi}{\mbox{\boldmath$\xi$}}
\newcommand{\brho}{\mbox{\boldmath$\rho$}}
\newcommand{\bbeta}{\mbox{\boldmath$\beta$}}
\newcommand{\bzeta}{\mbox{\boldmath$\zeta$}}

\def\bk{\boldsymbol{\kappa}}
\def\bmu{\boldsymbol\mu}
\def\bxi{\boldsymbol{\xi}}
\def\bz{\boldsymbol{\zeta}}

\def\ba{{\bf a}}
\def\bb{{\bf b}}
\def\bc{{\bf c}}
\def\be{{\bf e}}
\def\bff{{\bf f}}
\def\bg{{\bf g}}
\def\bn{{\bf n}}
\def\bp{{\bf p}}
\def\bq{{\bf q}}
\def\bs{{\bf s}}
\def\bt{{\bf t}}
\def\bu{{\bf u}}
\def\bv{{\bf v}}
\def\bw{{\bf w}}
\def\bx{{\bf x}}
\def\by{{\bf y}}
\def\bzz{{\bf z}}

\def\bD{{\bf D}}
\def\bE{{\bf E}}
\def\bF{{\bf F}}
\def\bH{{\bf H}}
\def\bJ{{\bf J}}
\def\bV{{\bf V}}
\def\bU{{\bf U}}
\def\bW{{\bf W}}
\def\bX{{\bf X}}
\def\bY{{\bf Y}}

\def\cA{{\cal A}}
\def\cC{{\cal C}}
\def\cD{{\cal D}}
\def\cE{{\cal E}}
\def\cF{{\cal F}}
\def\cG{{\cal G}}
\def\cI{{\cal I}}
\def\cJ{{\cal J}}
\def\cK{{\cal K}}
\def\cL{{\cal L}}
\def\cO{{\cal O}}
\def\cP{{\cal P}}
\def\cQ{{\cal Q}}
\def\cR{{\cal R}}
\def\cS{{\cal \Sigma}}
\def\cT{{\cal T}}
\def\cU{{\cal U}}
\def\cV{{\cal V}}

\def\scT{{_\cT}}
\def\sD{{_D}}
\def\sE{{_E}}
\def\sF{{_F}}
\def\sFz{{_{F_z}}}
\def\sK{{_K}}
\def\sI{{_I}}
\def\sb{{_b}}
\def\sN{{_N}}

\def\curl{{{\bf curl} \ }}
\def\rot{{\mbox{rot}\ }}
\def\BPI{{\bf \Pi}}

\def\cth{\cT_h}
\def\ctH{\cT_H}

\def\tJ{\tilde{\J}}

\def\hK{\widehat{K}}
\def\hx{\widehat{x}}
\def\hy{\widehat{y}}
\def\bhv{\widehat{\bv}}

\def\l{\ell}
\def\bl{\boldsymbol{\ell}}
\def\col{\colon}
\def\f12{\frac12}
\def\dfrac{\displaystyle\frac}
\def\dint{\displaystyle\int}
\def\nab{\nabla}
\def\p{\partial}
\def\sm{\setminus}
\def\dsum{\displaystyle\sum}
\newcommand{\pp}[2]{\frac{\partial {#1}}{\partial {#2}}}
\def\bzero{{\bf 0}}

\def\divv{\nab\cdot}
\def\divx{\nab_x\cdot}
\def\divtx{\nab_{t,x}\cdot}
\def\nabx{\nab_x}

\newcommand{\grad}{\nabla}
\newcommand{\curlt}{{\nabla \times}}
\newcommand{\gperp}{\nabla^{\perp}}
\newcommand{\gradt}{\nabla\cdot}

\def\forallqq{\quad\forall\,}
\def\aph{A^{1/2}}
\def\amh{A^{-1/2}}

\def\osc{{\rm osc \, }}

\def\Im{{\rm Im}}
\newcommand{\tr}{{\rm tr}}
\def\divvr{{\rm div}}
\def\curllr{{\rm curl}}
\def\curll{{\rm curl}}
\def\curl{{\bf curl}}
\newcommand{\bgrad}{{\bf grad}}
\newcommand\diam{\mathrm{diam\,}}
\renewcommand\Im{\mathrm{Im\,}}
\def\Span{\mbox{Span}}
\def\supp{\mbox{supp\,}}
\newcommand{\trace}{{\rm trace}}

\newcommand{\tri}{|\!|\!|}
\newcommand{\ljump}{\lbrack\!\lbrack}
\newcommand{\rjump}{\rbrack\!\rbrack}
\newcommand{\bdm}{\begin{displaymath}}
\newcommand{\edm}{\end{displaymath}}
\newcommand{\beq}{\begin{equation}}
\newcommand{\eeq}{\end{equation}}
\newcommand{\beqa}{\begin{eqnarray}}
\newcommand{\eeqa}{\end{eqnarray}}
\newcommand{\beqas}{\begin{eqnarray*}}
\newcommand{\eeqas}{\end{eqnarray*}}
\newcommand{\ul}{\underline}
\newcommand{\wh}{\widehat}
\newcommand{\la}{\langle}
\newcommand{\ra}{\rangle}

\newcommand{\Lt}{L^2(\Omega)}
\newcommand{\Lts}{L^2(\Omega)^2}
\newcommand{\Ltc}{L^2(\Omega)^3}
\newcommand{\Ho}{H^1(\Omega)}
\newcommand{\Hoh}{H^1(\wh{\Omega})}
\newcommand{\Hoi}{H^1(\Omega_i)}
\newcommand{\Hos}{H^1(\Omega)^2}
\newcommand{\Hoc}{H^1(\Omega)^3}
\newcommand{\Hoch}{H^1(\wh{\Omega})^3}
\newcommand{\Hoci}{H^1(\Omega_i)^3}
\newcommand{\Hoz}{H^1_0(\Omega)}
\newcommand{\Ht}{H^2(\Omega)}
\newcommand{\Hti}{H^2(\Omega_i)}
\newcommand{\Hts}{H^2(\Omega)^2}
\newcommand{\Htc}{H^2(\Omega)^3}
\newcommand{\Htz}{H^0(\Omega)}
\newcommand{\Hh}{H^{1/2}(\Gamma)}
\newcommand{\Hhi}{H^{1/2}(\Gamma_i)}
\newcommand{\Hmh}{H^{-1/2}(\Gamma)}
\newcommand{\Hdiv}{H(\divvr;\,\Omega)}
\newcommand{\Hdivh}{H(\divv;\,\wh \Omega)}
\newcommand{\hcurl}{H(\curl\,A;\,\Omega)}
\newcommand{\Hcurl}{H(\curll\,A;\,\Omega)}
\newcommand{\Hcrl}{H(\curll\,;\,\Omega)}
\newcommand{\hcrl}{H(\curl\,;\,\Omega)}
\newcommand{\Hcrlh}{H(\curll\,;\,\wh\Omega)}
\newcommand{\hcrlh}{H(\curl\,;\,\wh\Omega)}
\newcommand{\Wdiv}{\BW_0(\mbox{\divv}\,;\,\Omega)}
\newcommand{\Wcurl}{\BW_0(\mbox{\curl}\,A;\,\Omega)}
\newcommand{\WcrossV}{\BW \times V}

\begin{frontmatter}



\title{Battling Gibbs Phenomenon: \\
On Finite Element Approximations of Discontinuous Solutions of PDEs}


\author{Shun Zhang}
\ead{shun.zhang@cityu.edu.hk}
\address{Department of Mathematics, 
City University of Hong Kong, Hong Kong S.A.R., China}

\fntext[fn1]{This work was supported in part by Hong Kong Research Grants Council under the GRF Grant Projects No. CityU 11302519 and CityU 11305319.}



\begin{abstract}
In this paper, we want to clarify the Gibbs phenomenon when continuous and discontinuous finite elements are used to approximate discontinuous or nearly discontinuous PDE solutions from the approximation point of view.

For a simple step function, we explicitly compute its continuous and discontinuous piecewise constant or linear projections on discontinuity matched or non-matched meshes. For the simple discontinuity-aligned mesh case, piecewise discontinuous approximations are always good. For the general non-matched case, we explain that the piecewise discontinuous constant approximation combined with adaptive mesh refinements is a good choice to achieve accuracy without overshoots. For discontinuous piecewise linear approximations, non-trivial overshoots will be observed unless the mesh is matched with discontinuity. For continuous piecewise linear approximations, the computation is based on a "far-away assumption", and non-trivial overshoots will always be observed under regular meshes. We calculate the explicit overshoot values for several typical cases. Numerical tests are conducted for a singularly-perturbed reaction-diffusion equation and linear hyperbolic equations to verify our findings in the paper. Also, we discuss the $L^1$-minimization-based methods and do not recommend such methods due to their similar behavior to $L^2$-based methods and more complicated implementations. 
\end{abstract}

\end{frontmatter}

\section{Introduction}\label{intro}
\setcounter{equation}{0}

For a wide range of partial differential equations, their solutions can be discontinuous or nearly discontinuous. For example, for singularly perturbed problems, when the transient layer is very sharp, the solution can be viewed as discontinuous, see \cite{RST:08,FR:11}. For linear and nonlinear hyperbolic equations, discontinuous solutions can be caused by discontinuous initial data or shock forming, see \cite{HW:08,Hes:18}. For some cases, the location of discontinuity is known, for example, the discontinuity in the initial or boundary data. For many other cases, the exact location of discontinuity is unknown. When different numerical methods are used to solve such problems, the numerical solution is often {\em oscillatory} near a discontinuity, i.e., it overshoots (and undershoots). In spectral methods and Fourier analysis, it is called the Gibbs phenomenon, see \cite{HH:79,GS:97}. In some literature, they are also called {\em wiggles}, see \cite{Wesseling:01}. Such spurious oscillations or wiggles are unacceptable in many situations; see discussions on page 209 of Hesthaven \cite{Hes:18} for examples. These oscillations are essentially caused by a simple fact: continuous functions are used to approximate a discontinuous function. In numerical methods (finite difference, finite volume, spectral, and discontinuous Galerkin methods) for hyperbolic equations, some nonlinear tricks are used to eliminate or reduce the Gibbs phenomenon: artificial viscosity, limiters, filters, and other methods, see \cite{HW:08,Hes:18,GS:97}. 

On the other hand, in the numerical elliptic and singularly perturbed problems community using finite element methods, the Gibbs phenomenon and how to eliminate it are still quite mystical. Many numerical methods are tested for problems with discontinuity and near discontinuity. Sometimes the numerical solutions have no Gibbs phenomenon, while on the other occasions, the spurious oscillations appear for the same numerical method on a slightly different setting. There are many explanations for the phenomenon and how to reduce it. The following views are very common.

\begin{itemize}
\item [1] When the solution is discontinuous, if methods based on discontinuous piecewise polynomials such as discontinuous Galerkin (DG) methods are used, the overshoots can be eliminated or reduced. For example, in numerical computation for a singularly perturbed problem in a SIAM review paper \cite{FR:11} (p. 165), the authors claimed the DG methods would have no overshoots for almost discontinuous solutions. 

\item [2]
The adaptive finite element method with mesh refinements based on a posteriori error estimation is a useful procedure for detecting the regions and locations with bad approximations; see \cite{AO:00,Ver:13}. When the mesh is fine enough, the discontinuity is resolved, and the overshoot can be reduced. 

\item [3]
In the literature of finite difference/finite volume schemes, there are various {\em order barrier theorems}, see \cite{God:59,Wesseling:01}. For a simple linear equation ${d \phi}/{d x} = q(x)$,  Theorem 4.2.2 of \cite{Wesseling:01} states that linear discretization schemes of positive type are at most of first-order accurate. Similarly, there is the famous Godunov's order barrier theorem, see \cite{God:59} and Theorem 9.2.2 of \cite{Wesseling:01}: for the linear transport equation with constant coefficient, the linear one-step second-order accurate numerical schemes cannot be monotonically preserving.  These theorems are developed for specific equations and use the notions of positive or monotonically preserving numerical schemes. Theorem 4.2.2 of \cite{Wesseling:01} is proved by a simple local Taylor expansion and calculation of local truncation error for a one-dimensional uniform numerical grid. Do we have the same order barrier for other numerical methods and other partial differential equations?

\item [4]
There is also some research on the relation between the $L^1$-best approximation and the Gibbs phenomenon. The paper \cite{ST:99} of Saff and Tashevof shows that if continuous piecewise linear polynomials approximate a jump discontinuity on a discontinuity-matched uniform mesh in one dimension in the $L^1$-best approximation sense, the Gibbs phenomenon vanishes. In the context of numerical methods for PDEs, papers with the $L^1$-based numerical methods include \cite{Lavery:88,Lavery:89,Jiang:93,Guermond:04,HRZ:20,HRZ:19b,LD:20}.
\end{itemize}

In this paper, we want to clarify this question from the approximation point of view by using a simple model problem. Since most finite element or DG methods are based on projections or pseudo-projections in Hilbert spaces, the best result a numerical method can achieve of a certain approximation often cannot be better than its $L^2$-projection. That is, if $V_h$ is the discrete approximation space, the best we can hope for the numerical error is often in the following form (note that for many complicated problems, it is essential to use adaptive methods to achieve the desired optimality; and many methods cannot have such a good result, for example, the DG method for the linear hyperbolic equation):
\beq
\|u - u_h\|_0 \leq C \inf_{v_h \in V_h} \|u-v_h\|_0 = C \|u-\Pi_{V_h}u\|_0,
\eeq
where $\Pi_{V_h}$ is the $L^2$-projection operator onto the discrete approximation space $V_h$ and $\|\cdot\|_0$ is the $L^2$-norm. 
{Since the Gibbs phenomenon is about the wiggles/overshoots/undershoots, the error is pointwise. Thus, to study the Gibbs phenomenon, we should measure the error in an $L^\infty$ like norm. In Theorem 6.5 of \cite{CK:20},  for the singularly perturbed reaction-diffusion problem with a two-step finite element method, the following result, 
\beq \label{CK}
\|u -u_h\|_{L^\infty(K)} \leq  C\|u -\Pi_{V_h} u\|_{L^\infty(K)} + \mbox{extra terms on slightly larger domains of } K,
\eeq
is proved, where $V_h$ is the discontinuous piecewise polynomial space on a finite element mesh, $K$ is an element of the finite element mesh, and $\Pi_{V_h}$ is the $L^2$-projection onto $V_h$.  Inspired by \eqref{CK}, at least from the approximation point of view, to study the Gibbs phenomenon of finite element approximations to (nearly) discontinuous solutions, we need to study the term $\|u -\Pi_{V_h} u_h\|_{L^\infty(K)}$ for different $V_h$ spaces (continuous or discontinuous finite elements on discontinuity-aligned or nonaligned, uniform or adaptive meshes) with $\Pi_{V_h}$ as the $L^2$-projection.

We also need to emphasize that the numerical stability of the numerical method for a specific equation is different from the Gibbs phenomenon or spurious oscillation, as discussed on page 119 of Wesseling \cite{Wesseling:01}. The numerical stability is essential for the success of the numerical method, but it cannot guarantee the disappearance of the Gibbs phenomenon. The Gibbs phenomenon (spurious oscillation or wiggle) is essentially an approximation problem. Some tricks to increase the numerical stability, such as the artificial/numerical viscosity, will enhance both the stability and reduce the Gibbs phenomenon. We assume that the artificial/numerical viscosity is not used.

Choosing $u$ to be the simplest discontinuous function, a step function, we try to study the properties of $\Pi_{V_h}u$ with $V_h$ being the continuous or discontinuous finite element space on a discontinuity-aligned or nonaligned, uniform or adaptive mesh.  
Based on explicit computations, we find the following results. For the simple discontinuity-aligned mesh case, piecewise discontinuous approximations are always good in terms of overshoot. A special but important discontinuity-aligned mesh case is the boundary layer, where the (almost) discontinuity is on the boundary and is matched with the numerical mesh. For the general discontinuity-nonaligned mesh case, we explain that {\bf the piecewise discontinuous constant approximation combined with adaptive mesh refinements is the best choice to achieve accuracy without overshoot}. For discontinuous piecewise linear approximations, non-trivial overshoots will be observed unless the mesh is matched with discontinuity. For continuous piecewise linear approximations, non-trivial overshoots will always be observed under regular meshes. We calculate the explicit overshoot values for several typical cases.

We also find that a special high ratio mesh obtained by adaptive mesh refinements combined with coarsening can reduce the overshoot. However, such a case is not always possible in two and three dimensions or one dimension with other error sources.



The paper is organized as follows. Section 2 discusses a simple case of piecewise discontinuous polynomial approximations on a discontinuity-aligned mesh.  Section 3 describes a model problem of approximating a step function by $L^2$ projections onto piecewise constant and linear function spaces. Discontinuous piecewise constant and linear approximations are explicitly computed for the model problem. We also compute continuous piecewise linear approximations under a ”far-away assumption” for the model problem. In Section 4, linear conforming, linear DG, and lowest-order mixed finite element methods are tested on a one-dimensional singularly perturbed problem. Several test problems of 2D linear transport problems with P0- and P1-DGFEMs are done in Section 5. In Section 6, we comment on known numerical test results in the literature. We discuss the $L^1$-minimization-based methods in Section 7. We present a table to summarize the results for different cases in Section 8. In Section 9, we make some concluding remarks.

\section{Piecewise Discontinuous Polynomial Approximations on A Discontinuity-Aligned Mesh}
In this section, we consider an almost trivial but important case, piecewise discontinuous polynomial approximations on a discontinuity-aligned mesh. Assume $\cT_h =\{K\}$ is a finite element mesh on $d$-dimensional polygonal/polyhedral domain, and the function $u|_K \in H^{2}(K)$ on each element $K\in\cT_h$. By the Sobolev embedding theorem, we gave $u|_K \in C^{0}(K)$. The norm $\|u\|_{L^{\infty}(K)}$ is also well-defined. The function is continuous inside each element, and thus the discontinuity of the function (if any) is aligned with the mesh.  Let $P_k(K)$ be the space of polynomials defined on $K$ whose degree is less or equal to integer $k \geq 0$ and let $P_k(\cT_h):=\{ v \in L^2(\O) : v|_K = P_k(K), \forall K\in\cT_h\}$ be the space of piecewise discontinuous polynomial space. The diameter of the element $K$ is defined as $h_K$. The notation $|\cdot|_{2,K}$ donotes the $H^2$-semi-norm on an element $K$.
 We have the following Theorem.

\begin{thm}
Assume on each element $K\in \cT_h$, $u|_K \in H^{2}(K)$. Let $u_{k,K}\in P_k(K)$ be the $L^2$-projection of $u$ in $P_k(K)$, then for some constant $C>0$ independent of $h$,
\beq
\|u-u_{k,K}\|_{L^{\infty}(K)} \leq C h_K^{2-d/2} |u|_{2,K} \quad \mbox{for  } d= 1,2,3.
\eeq
\end{thm}
\begin{proof}
The proof here is a modification of the argument in p.93 of Braess \cite{Braess:07}.
Since $u\in C^0(K)$, let $u_{I,K}$ be the well-defined nodal interpolation of $u$ in $P_k(K)$. By the triangle inequality, we have
\beq \label{uukk}
\|u-u_{k,K}\|_{L^{\infty}(K)} \leq \|u-u_{I,K}\|_{L^{\infty}(K)} + \|u_{I,K}-u_{k,K}\|_{L^{\infty}(K)}. 
\eeq
By the inverse estimate (Lemma 1.138 of \cite{EG:04}), we get
\beq
\|u_{I,K}-u_{k,K}\|_{L^{\infty}(K)}\leq Ch_K^{-2/d}\|u_{I,K}-u_{k,K}\|_{0,K}.
\eeq
By the triangle inequality and the approximation properties of the $L^2$-projection and the nodal interpolation, 
\beq
\|u_{I,K}-u_{k,K}\|_{0,K} \leq \|u-u_{k,K}\|_{0,K}+\|u-u_{I,K}\|_{0,K} \leq C h_K^{2}|u|_{2,K}.
\eeq
Combining above two inequalities, we have
\beq
\|u_{I,K}-u_{k,K}\|_{L^{\infty}(K)} \leq C h_K^{2-2/d}|u|_{2,K}.
\eeq
The theorem is proved by substituting it into \eqref{uukk} and using the approximation property of the interpolation (Theorem 16.1 of \cite{Ciarlet:91}):
$$
\|u-u_{I,K}\|_{L^{\infty}(K)} \leq C h_K^{2-2/d}|u|_{2,K}.
$$
\end{proof}

\begin{rem}
Note that there are various more refined max-norm/pointwise error estimates for different methods and equations; see \cite{Ciarlet:91,Wahlbin:91,JP:86}. Most of these analyses are based on the global smoothness of the solution.  
\end{rem}

\begin{rem}
From the theorem, it is clear that for the $L^2$-projection of the piecewise discontinuous polynomial approximation on a discontinuity-aligned mesh, the pointwise error will decrease to zero as the mesh size decreases. Combined with adaptive mesh refinements, there is no or only neglectable Gibbs phenomenon for this case.
\end{rem}

\begin{rem}
Note that some boundary layer problems fall into this category. When there is a sharp solution change near the boundary, the solution can be viewed as discontinuous on the Dirichlet boundary. In this case, the mesh is always aligned with the boundary (assuming the domain is polygonal/polyhedral); thus it is a discontinuity-aligned case.

From the discussion in this section, to eliminate the Gibbs phenomenon/overshoot, the boundary condition should be enforced weakly. The famous Nitsche's method \cite{Nitsche:71} and discontinuous Galerkin method \cite{ABCM:02,HW:08} are successful for this case. We also need to emphasize that weak enforcement should be chosen carefully. A simple penalization may not work since it is not weak enough, see some numerical experiments in \cite{Lin:08}.
\end{rem}

\begin{rem}
For the Gibbs phenomenon near a discontinuity, we are only interested in the value of the numerical solution that overshoots/undershoots. It may be possible that the numerical approximation lies between the discontinuity gap and the maximum error at that location is large, but no Gibbs phenomena happen. Thus, the measurement we use to quantify the Gibbs phenomena for a nonaligned mesh is a little weaker than the standard $L^\infty$-norm; see also our definition of $\mathtt{os}$ in Sections 3.2.2 and 3.3.2.
\end{rem}

\section{An Illustrative One-Dimensional Problem}
\setcounter{equation}{0}

\subsection{A Model Problem}
Consider a step function $u$: 
\beq \label{u}
u(x) = \left\{\begin{array}{cc}
- 1 & x<0, \\[1mm]
1 & x>0.
	\end{array}\right. 
\eeq
Here, the discontinuity gap $u(0+) - u(0-)$ is $2$. For more general cases with a gap $2c$, $c>0$, it is easy to see that the corresponding results below are proportional to $c$.

Let $\cT_h =\{K\}$ be a one-dimensional mesh on $\O = (-a,b)$ with elements (intervals in this 1-D setting) denoted by $K$. We assume that $a$ and $b$ are big enough positive numbers with respect to the mesh size. 
Consider three approximation spaces defined on $\cT_h$:
\begin{eqnarray}
P_0(\cT_h) &:=& \{ v \in L^2(\O) : v|_K = P_0(K), \forall K\in\cT_h\},\\
P_1(\cT_h) &:=& \{ v \in L^2(\O) : v|_K = P_1(K), \forall K\in\cT_h\},\\
S_1(\cT_h) &:=& \{ v \in C^0(\O) : v|_K = P_1(K), \forall K\in\cT_h\}.
\end{eqnarray}
Let $V_h$ be $P_0(\cT_h)$, $P_1(\cT_h)$, or $S_1(\cT_h)$. The $L^2$-projection $u_h\in V_h$ is defined by:
$$
(u_h, v_h)_{\O} = (u,v_h)_{\O}\quad \forall v_h \in V_h.
$$
We use notations $u_{h,0}$, $u_{h,1}$, and $u_{c,1}$ to denote the $L^2$-projections of $u$ on $P_0(\cT_h)$, $P_1(\cT_h)$, and $S_1(\cT_h)$, respectively.

\subsection{Discontinuous Approximations}
We consider two discontinuous approximations first: $P_0$ and $P_1$ approximations. Note that, for a discontinuous projection, the approximation is exact on those elements away from the discontinuity. Thus, we only need to discuss the special interval that contains the discontinuity. Let $0\leq t \leq 1$ and $h>0$. We consider the interval $I_0 = (-th, (1-t)h)$, see (a) of Fig.  \ref{P1_non}.

\subsubsection{Discontinuous Piecewise Constant Approximation}
A simple computation shows that  the constant projection of $u$ on the interval $I_0$ is the following,
\beq
	u_{h,0}|_{I_0} = 1-2t.
\eeq
It is obvious that $ -1 \leq 1-2t\leq 1$ for $t\in [0,1]$, thus, there is no overshoot. If $t=0$ or $1$, which is a discontinuity-aligned mesh, the numerical approximation is exact.

Computing the $L^2$-error in the interval $I_0$, we have,
\beq
\|u-u_{h,0}\|_{0, I_0} = 2 \sqrt{t(1-t)h}.
\eeq
The worst case is $\|u-u_{h,0}\|_{0, I_0} = \sqrt{h}$, when $t=1/2$. This also matches the a priori analysis that $\|u-u_{h,0}\|_{0,I_0} \leq C h^{1/2-\epsilon}$ for $u\in H^{1/2-\epsilon}(I_0)$ for an arbitrary small $\epsilon >0$.

With a good a posteriori error estimator that can identify the bad approximated elements, the discontinuity-crossing elements will be found and divided, and the  $L^2$- and other integration-based norms of error will become smaller and smaller. 

\subsubsection{Discontinuous Piecewise Linear Approximation}\label{S322}
In this subsection, we consider the $L^2$-projection of $u$ onto a linear function on $I_0$. To this end, let $\lambda_{-1}$ and $\lambda_1$ be the linear Lagrange basis functions define on $(-th,(1-t)h)$ with $\lambda_{i}(z_i)=1$ and $\lambda_{i}(z_{-i})=0$, where $i= 1$ or $-1$, and $z_{-1} = -th$ and $z_1 = (1-t)h$.

Let the projection $u_{h,1}|_{I_0} = U_{-1}\lambda_{-1} + U_{1}\lambda_{1}$. The projection problem is 
$$
\left( \begin{array}{cc}
(\lambda_{-1}, \lambda_{-1})_{I_0} & (\lambda_{1}, \lambda_{-1})_{I_0}\\[2mm]
(\lambda_{-1}, \lambda_{1})_{I_0} & (\lambda_{1}, \lambda_{1})_{I_0}
\end{array}
 \right)
 \left( \begin{array}{c}
U_{-1}\\[2mm]
U_1
\end{array}
 \right) =
\left( \begin{array}{c}
(u,\lambda_{-1})_{I_0}\\[2mm]
(u,\lambda_{1})_{I_0}
\end{array}
\right).
$$
All terms of the matrix problem can be computed exactly. Solving the projection problem, we get
\begin{eqnarray*}
U_{-1} = 1-8t+6t^2\quad \mbox{and}\quad U_1 =  1+4t-6t^2.
\end{eqnarray*}
It is easy to show that
\begin{eqnarray*}
(U_{1}-1)  > -(U_{-1}+1) \geq 0 && \mbox{  if  } 0<t<1/2, \\[1mm]
 -(U_{-1}+1) >(U_{1}-1)  \geq 0 && \mbox{  if  } 1/2<t<1.
\end{eqnarray*}
To quantify the overshoot phenomena, we define the following overshoot value function:
$$
\mathtt{os} = \max(U_{1}-1, -(U_{-1}+1), 0).
$$
As discussed in Remark 2.5, we do not use the  $\|u-u_{h,1}\|_{L^{\infty}}$ here since we are only  interested in the $\mathtt{os}$ (which is weaker than  $\|u-u_{h,1}\|_{L^{\infty}}$).

We plot the value of $\mathtt{os}$ for $0\leq t \leq 1$ on the right of Fig. \ref{P1_non}. It is easy to see that only when the mesh is aligned with the discontinuity ($t=0$ or $1$), $\mathtt{os}$ is zero. For $t$ away from $0$ or $1$, the overshoot phenomenon is severe. The maximum overshoot value $2/3$  appears at $t=1/3$ or $2/3$.

We also notice that if the bisection of mesh is used, the relative position of the discontinuity normally will not converge to the left or the right point of the internal. For example, consider the interval $(-1/3,2/3)$, i.e., the discontinuity is at the $1/3$ position of the internal. If we bisect the interval, the new interval containing the discontinuity is $(-1/3, 1/6)$, and the discontinuity is at the $2/3$ position. Keeping doing the bisection, we will find that the discontinuity will jump between the $1/3$ and $2/3$ positions. In this case, $t$ will always be $1/3$ or $2/3$, and the overshoot value will not decrease. For other initial positions other than the aligned case, the overshoot values will oscillate between $[\delta, 2/3]$, for some $0<\delta<2/3$. 

\begin{figure}[!htb]
\centering 
\subfigure[the interval $I_0$]{ 
\includegraphics[width=0.45\linewidth]{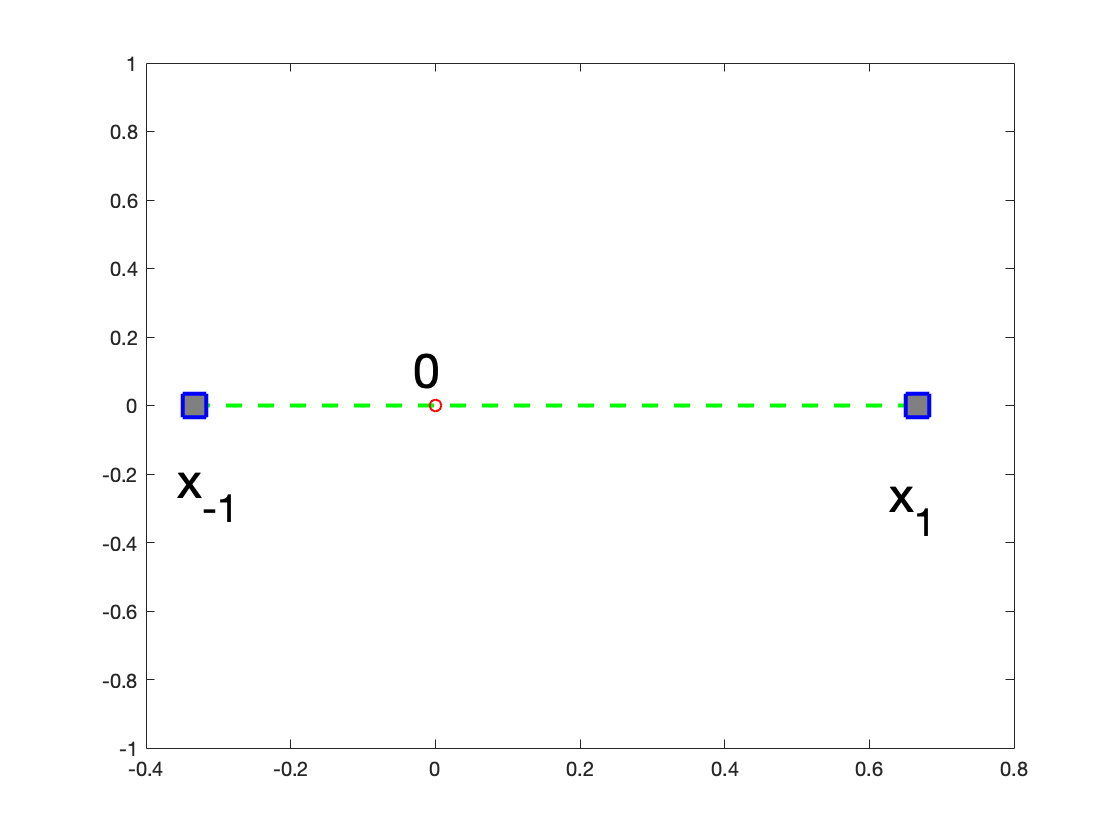}}
~
\subfigure[overshoot values as $t \in (0,1)$]{
\includegraphics[width=0.45\linewidth]{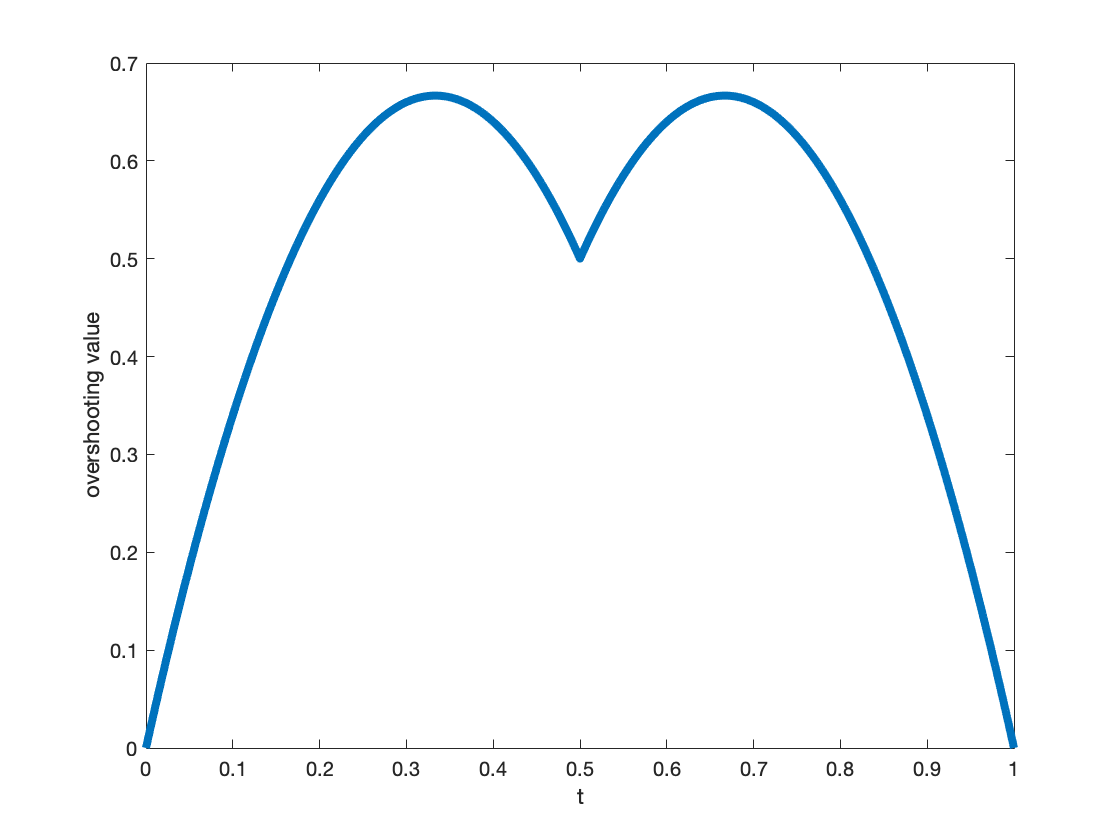}}
\caption{Discontinuous piecewise linear approximation on a discontinuity non-aligned  mesh for a step function}
 \label{P1_non}
\end{figure}
%
Computing the $L^2$-error in the interval $I_0$, we have,
\beq
\|u-u_{h,1}\|_{0, I_0} = 2\sqrt{t(1-t)(1-3t+3t^2)h} .
\eeq
When $t=0$ or $1$, the error is zero. The largest error $\sqrt{h/3}$ appears at $t= 0.5\pm \sqrt{3}/6$.
Note that $1/4\leq 1-3t+3t^2 \leq 1$ for $0\leq t\leq 1$, we do have 
\beq
\|u-u_{h,1}\|_{0,I_0} \leq \|u-u_{h,0}\|_{0, I_0}. 
\eeq
Thus, for discontinuous approximations, even though the $L^2$-norm of error of the $P_1$ approximation is better than the $P_0$ approximation, they are of the same order if the mesh is not aligned, and the overshoot cannot be avoided in most cases.

\subsection{Continuous Piecewise Linear Approximation}
For the continuous piecewise linear approximation, the exact solution is not in the approximation space, even for the discontinuity-aligned mesh case. The projection is not local to one element anymore. Luckily, based on the numerical experiments, we found that the Gibbs phenomenon is most severe on the two nodes around the discontinuity point. The more nodes between the discontinuity and the node to be observed, the smaller the error. We assume the following {\bf far-away assumption} for the continuous piecewise linear approximation.

\begin{asm}
{\bf Approximation on nodes far-away from the discontinuity (far-away assumption):} For a node $x_k$ of the mesh $\cT_h$, if there is at least two nodes between the discontinuity position and $x_k$, the error at $x_k$, $|u_{c,1}(x_k) - u(x_k)|$, is assumed to be less than 1/10 the maximum error at the possible overshoot locations (the discontinuity point, and its two neighbor nodes). 

For a node $x_k$ of the mesh $\cT_h$, if there is one node between the discontinuity position and $x_k$, the error at $x_k$, $|u_{c,1}(x_k) - u(x_k)|$, is assumed to be less than 1/2 of  the maximum error at the possible overshoot locations (the discontinuity point, and its two neighbor nodes).
\end{asm}
In our calculations, we plan to study the maximum error at the possible overshoot locations (the discontinuity point, and its two neighbor nodes) in our calculations. We assume that the error at $x_k$, which is two nodes away from the discontinuity, is neglectable. For those points only one node away from the discontinuity, we will either compute their values or assume it has a perturbation $\epsilon$ to be determined later.

\subsubsection{Discontinuity-aligned mesh}
We first consider a possible ideal case: the mesh is aligned with the discontinuity and the mesh is symmetric with respect to the discontinuity at $0$. 
Let $\cT_h$ be a symmetric mesh on $[-1,1]$ with nodes: $x_{-N} < x_{-N+1} <\cdots< x_{-1} < x_0 < x_1<\cdots<x_{N-1}<x_N$, where $x_{-i} = -x_i$ and $x_0=0$. 

\begin{figure}[!htb]
\centering 
\subfigure[$h=2^{-4}$]{ 
\includegraphics[width=0.45\linewidth]{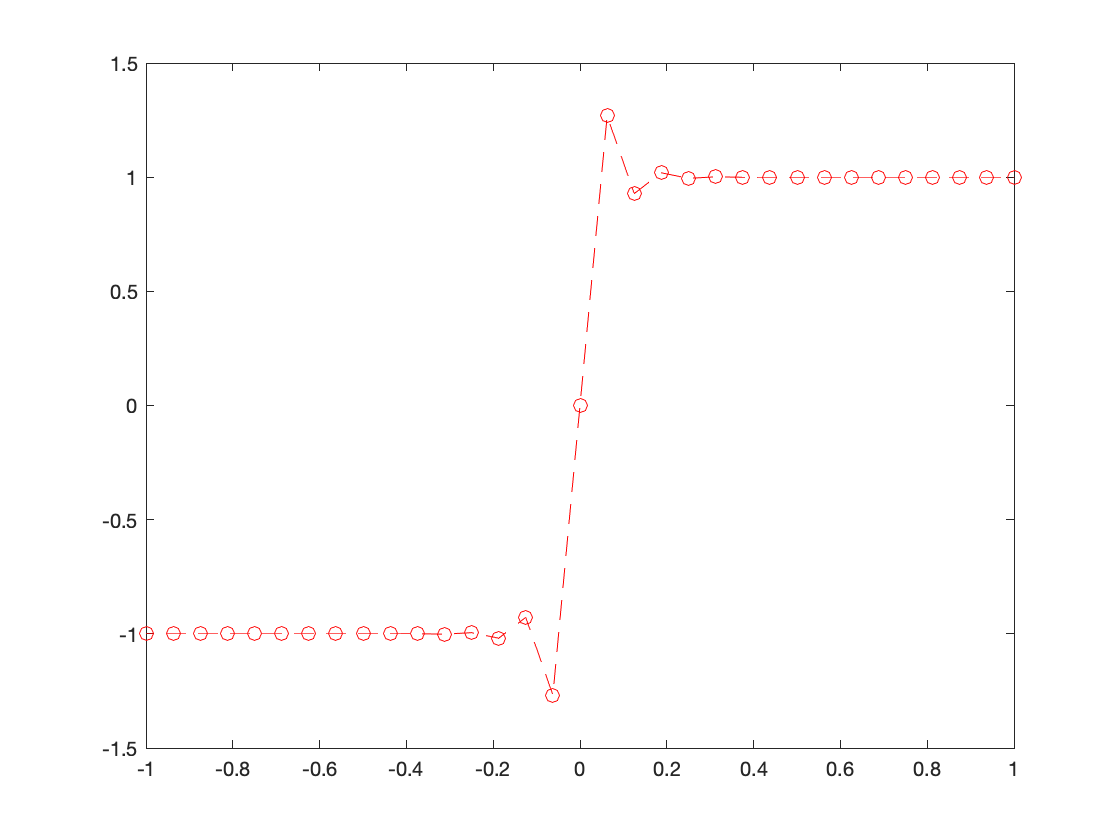}}
~
\subfigure[$h=2^{-14}$]{
\includegraphics[width=0.45\linewidth]{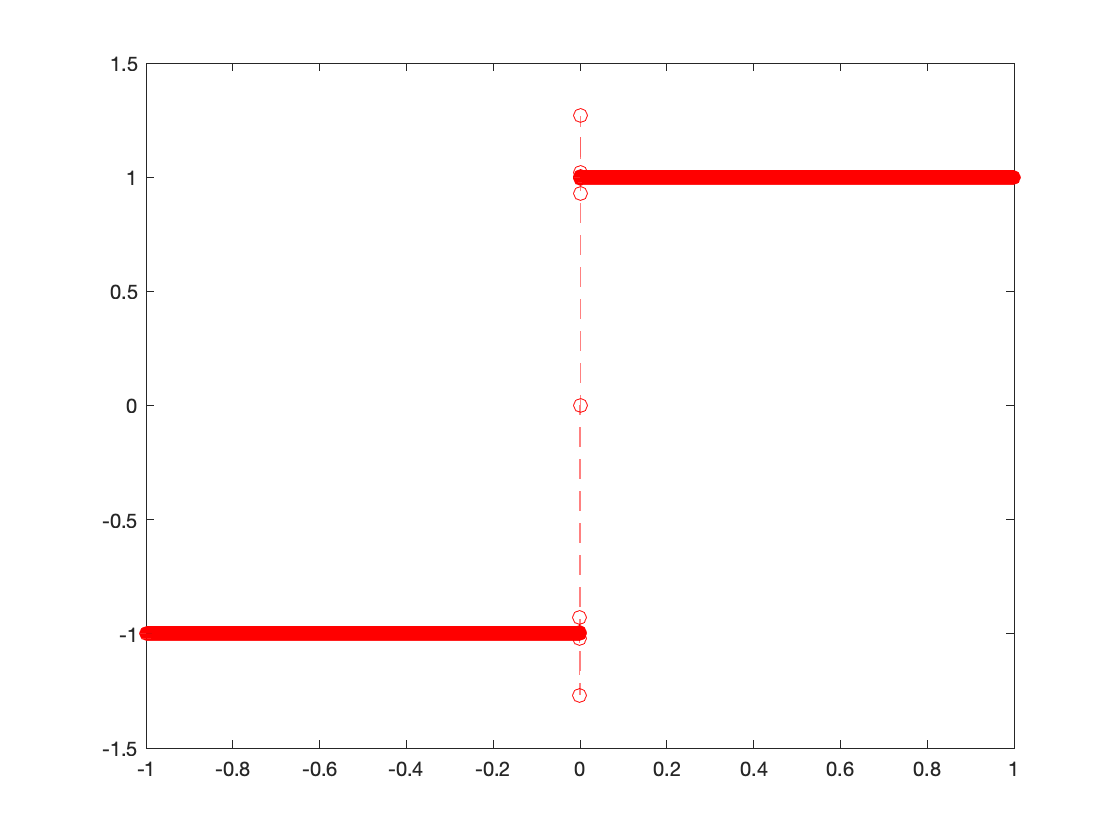}}
\caption{Continuous piecewise linear approximations on a uniformly refined mesh for a step function}
 \label{CP1_uniform}
\end{figure}

We first test on a uniform mesh with $h=1/N$ and $x_{\pm i} = \pm i h$. From the numerical tests, for example,  $h=1/16$ and $h=2^{-14}$ in Fig. \ref{CP1_uniform}, we note that the values of $u_{c,1}$ at the center nodes $x_{\pm i}$ are identical for different choice $h$ for $i \leq 4$. For $i$ big enough, $u_{c,1}(x_{\pm i})$ is almost identical to $u$. The values $u_{c,1}(x_{i})$ for $i= -4,-3, \cdots, 3, 4$ are $-0.9948$, $-1.0192$, $-0.9282$, $-1.2679$, $0$, $1.2679$, $0.9282$,   $1.0192$, and $0.9948$. The severe overshoots only appear on the two nodes around the discontinuity with a size of $0.2679$. For the points $i=-3$ and $i=3$, the errors are $0.0192$, less than $1/10$ of the overshoot error of $0.2679$. For the points $i=-2$ and $i=2$, the errors are $0.0718$, less than $1/3$ of the overshoot error $0.2679$. The far-away assumption is true in this case.

With the far-away assumption, a simple calculation can be done to show why the values of $u_{c,1}$ at $x_{\pm 1}$ are close to $\pm 1.25$ as observed. We simplify the calculation by assuming that $u_{c,1} =u$ at all $x_{\pm i}$ for $i> 2$.  We assume the values of $u_{c,1}(x_{\pm 2}) = \pm (1+\epsilon)$, where $|\epsilon|$ is a small number compared to $1$.  We have the rest three values at $x_{\pm 1}$ and $x_0$ to be computed by the projection. By the symmetry of nodes, it is easy to see that $u_{c,1}$ at $x_0 = 0$ is $0$ and $u_{c,1}(x_i) = - u_{c,1}(x_{-i})$.   Thus we only need to determine the value of $u_{c,1}(x_1) =: U_1$.  Let $\lambda_1$ and $\lambda_2$ be the linear Lagrange basis functions on the mesh, with $\lambda_{i}(x_j) = \delta_{i,j}$, $i,j \in \{1,2\}$. Then $u_{c,1}$ on $(x_0,x_2)$ can be written as 
$$
u_{c,1} = U_1 \lambda_1 + (1+\epsilon) \lambda_2, \quad x\in(x_0,x_2),
$$
The value of $U_1$ can be obtained by a simple projection,
$$
( U_1 \lambda_1 + (1+\epsilon) \lambda_2,\lambda_1)_{(x_0,x_2)} = (1, \lambda_1)_{(x_0,x_2)}.
$$   
A simple calculation shows
\beq
U_1 = 1.25 -\epsilon/4.
\eeq
Thus, the overshoot value is always close to $1/8$ of the discontinuity jump (2 in our example).  The exact value of $1.2679$ is due to the fact that $\epsilon$ is about $-0.0718$ in the example.

The $L^2$-errors in the two adjacent elements of $0$ are also easy to be computed by choosing $\epsilon$ to be $0$,
\beq \label{L2matchCP1}
\|u-u_{c,1}\|_{0, (0,h)} \approx \|u-u_{c,1}\|_{0, (-h,0)} \approx \sqrt{13h/48} \approx  0.5204 \sqrt{h}. 
\eeq

\subsubsection{Discontinuity-nonaligned mesh}
We consider the case where the mesh is not aligned with the discontinuity. From the above discussion, we notice that the size of each element is not essential (their ratio is important, see the example below). Three cases are considered: a local uniform mesh, a local adaptive mesh, and an adaptive mesh with coarsening.

\paragraph {\bf A local uniform mesh}
Suppose $0<t<1$ and $h>0$, we consider the following mesh of 5 elements with size $h$: $x_{-3} = -(t+2)h$, $x_{-2} = -(t+1)h$, $x_{-1} = -th$, $x_{1}=(1-t)h$, $x_2=(2-t)h$, and $x_3 = (3-t)h$, see the left of Fig. \ref{CP1_non} for a test mesh with $t=1/3$ and $h=1$. The discontinuity cuts through the central interval $I_0 = (-th, (1-t)h)$ at $0$. Let $\lambda_{i}$ be the linear Lagrange basis function for $i =\pm 1,\pm 2, \pm 3$. As before, we assume that $u_{c,1}$ is exact at $i = \pm 3$ (we omit the small perturbation $\epsilon$ for simplicity). Then on the interval $(x_{-3},x_3)$, 
$$
u_{c,1} = -\lambda_{-3} + U_{-2} \lambda_{-2}  + U_{-1} \lambda_{-1}+ U_{1} \lambda_{1} + U_{2} \lambda_{2}+\lambda_{3}.
$$ 
Solving the following projection problem, 
$$
(u_{c,1}, \lambda_{i})_{(x_{-3},x_3)} =  (u, \lambda_{i})_{(x_{-3},x_3)}, \quad i= \pm 1, \pm 2,
$$
we get
\begin{eqnarray*}
U_{-2} &=& -3 (87 - 60 t + 38 t^2)/209, \quad
U_{-1} = (-1 - 720 t + 456 t^2)/209,\\
U_1 &=&  (265 + 192 t - 456 t^2)/209,\quad
U_2 = 3 (65 - 16 t + 38 t^2)/209.
\end{eqnarray*}
We define the following overshoot value function:
$$
\mathtt{os} = \max(-(U_{-2}+1), -(U_{-1}+1),(U_{1}-1),(U_{2}-1),0).
$$
We plot the overshoot value with respect to $t$ on the right of Fig.\ref{CP1_non}. The overshoot values lie in the interval 
$$[0.1818, 0.3646].$$ 
Also, when $0<t<1/2$, the maximum overshoot appears at $U_1$, and when $1/2<t<1$, the maximum overshoot happens at $U_{-1}$, that is, the maximum overshoot appears at one of the endpoints of $I_0$, which is farther from the discontinuity. The overshoot is non-trivial for any $0<t<1$.

\begin{figure}[!htb]
\centering 
\subfigure[a test mesh with $t=1/3$ and $h=1$]{ 
\includegraphics[width=0.45\linewidth]{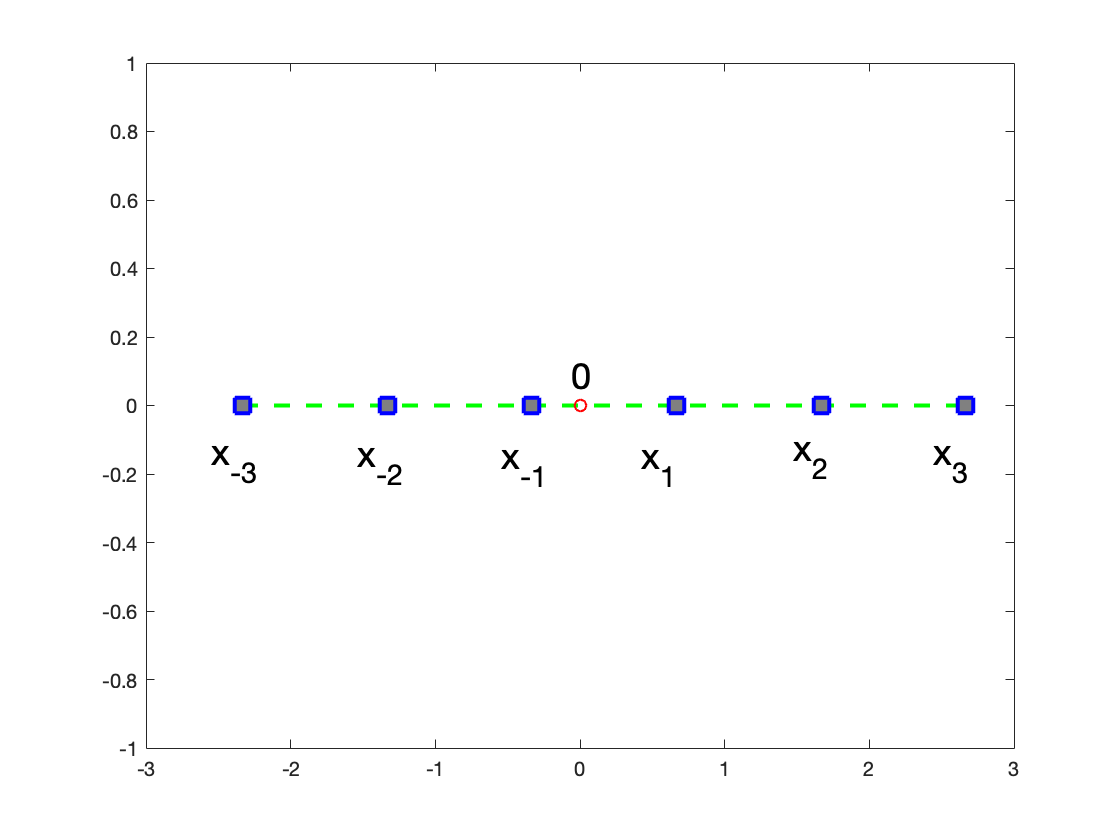}}
~
\subfigure[overshoot values as $t \in (0,1)$]{
\includegraphics[width=0.45\linewidth]{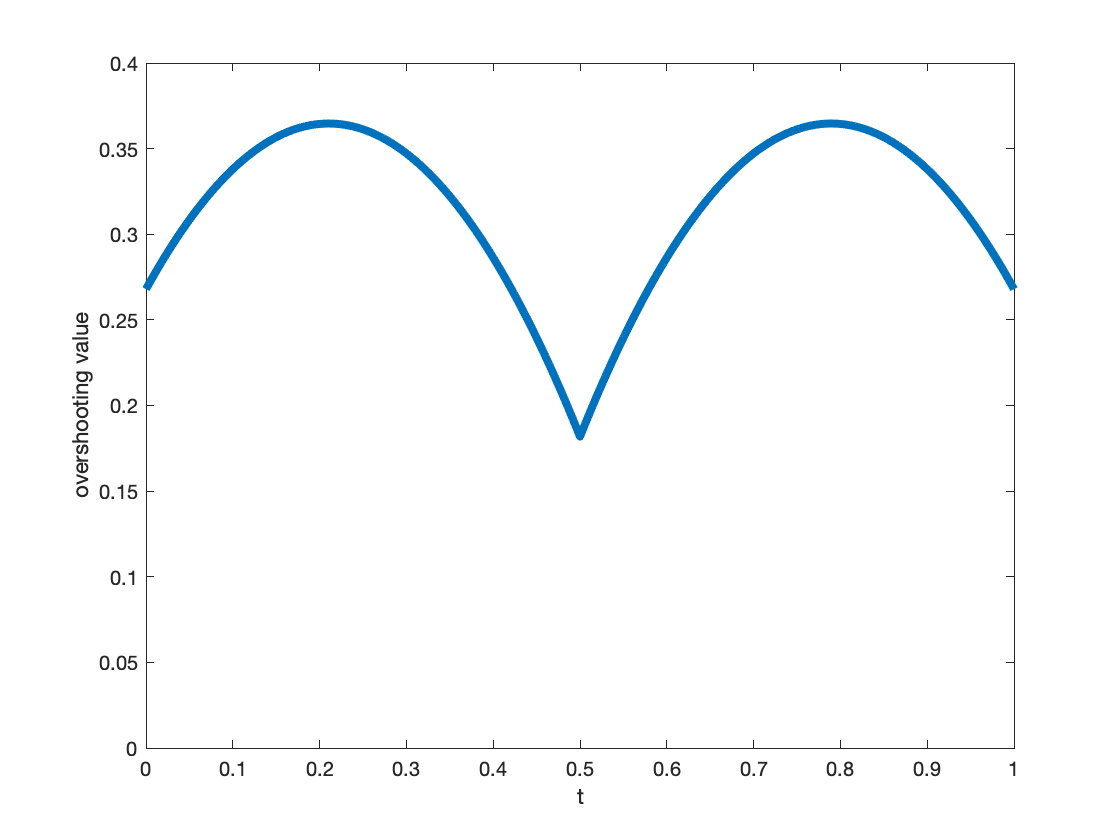}}
\caption{Continuous piecewise linear approximation on a discontinuity-nonaligned  mesh for a step function}
 \label{CP1_non}
\end{figure}
We have the following $L^2$-error:
\beq
\|u-u_{c,1}\|_{0, I_0} = 2\sqrt{h(8869 + 60189 t - 294117 t^2 + 467856 t^3 - 233928 t^4)/131043} .
\eeq
For $0\leq t \leq 1$, we get $\|u-u_{c,1}\|_{0, I_0}$ is between $0.5203 \sqrt{h}$ and $0.6236 \sqrt{h}$. Note this error matches with the result of \eqref{L2matchCP1}, which corresponds to $t=0$ or $1$.

\paragraph{\bf A local adaptive mesh}
We also test the following mesh: $x_{-3} = -(t+4)h$, $x_{-2} = -(t+2)h$, $x_{-1} = -t h$, $x_{1}=(1-t)h$, $x_2=(2-t)h$, and $x_3 = (4-t)h$, with sizes $2h$, $h$, $h$, and $2h$, respectively. This kind of mesh can be viewed as an example of an adaptively refined mesh. 
With this mesh setting and assuming the values of $u_{c,1}$ at $x_{\pm 3}$ are exact,  by a similar procedure, we get
\begin{eqnarray*}
U_{-2} &=& 1/241 (-281 + 138 t - 87 t^2),\\
U_{-1} &=& 3/241 (-27 - 184 t + 116 t^2),\\
U_1 &=& 1/241 (325 + 144 t - 468 t^2),\\
U_2 &=& 1/241 (227 - 24 t + 78 t^2).
\end{eqnarray*}
We show the result on the right of Fig.\ref{CP1_non2}. The overshoot value lies in the interval 
\beq \label{os_cp1}
[0.1358, 0.3954].
\eeq 
Also, as before, we find that the maximum overshoot happens at $U_1$ when $0<t<1/2$,  and  it happens at $U_{-1}$ when $1/2<t<1$. For any choice of $0<t<1$, the overshoot is non-trivial. We vary the element size of elements adjacent to $I_0$ from $h$, $2h$, $4h$, to $8h$, and get similar results of non-trivial overshoots. This shows that adaptive mesh refinements by bisection will not make the overshoot phenomenon disappear for a continuous linear approximation.

\begin{figure}[!htb]
\centering 
\subfigure[a test mesh]{ 
\includegraphics[width=0.45\linewidth]{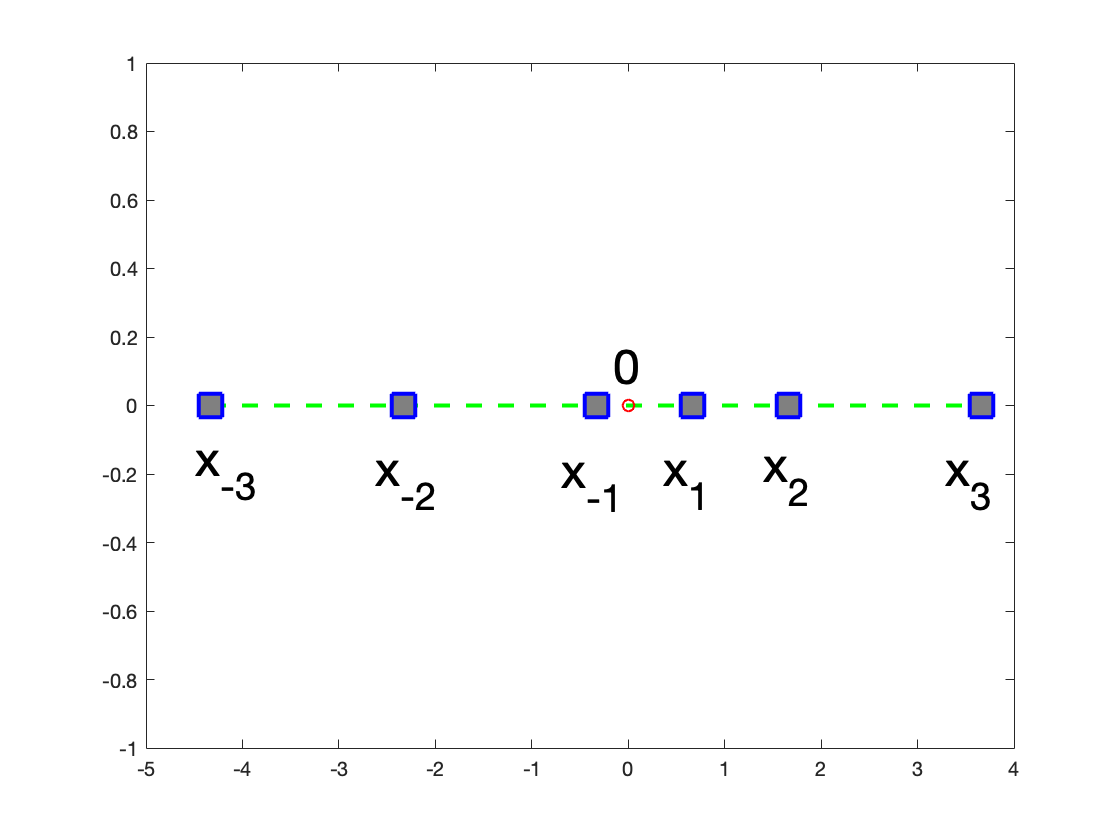}}
~
\subfigure[overshoot values as $t \in (0,1)$]{
\includegraphics[width=0.45\linewidth]{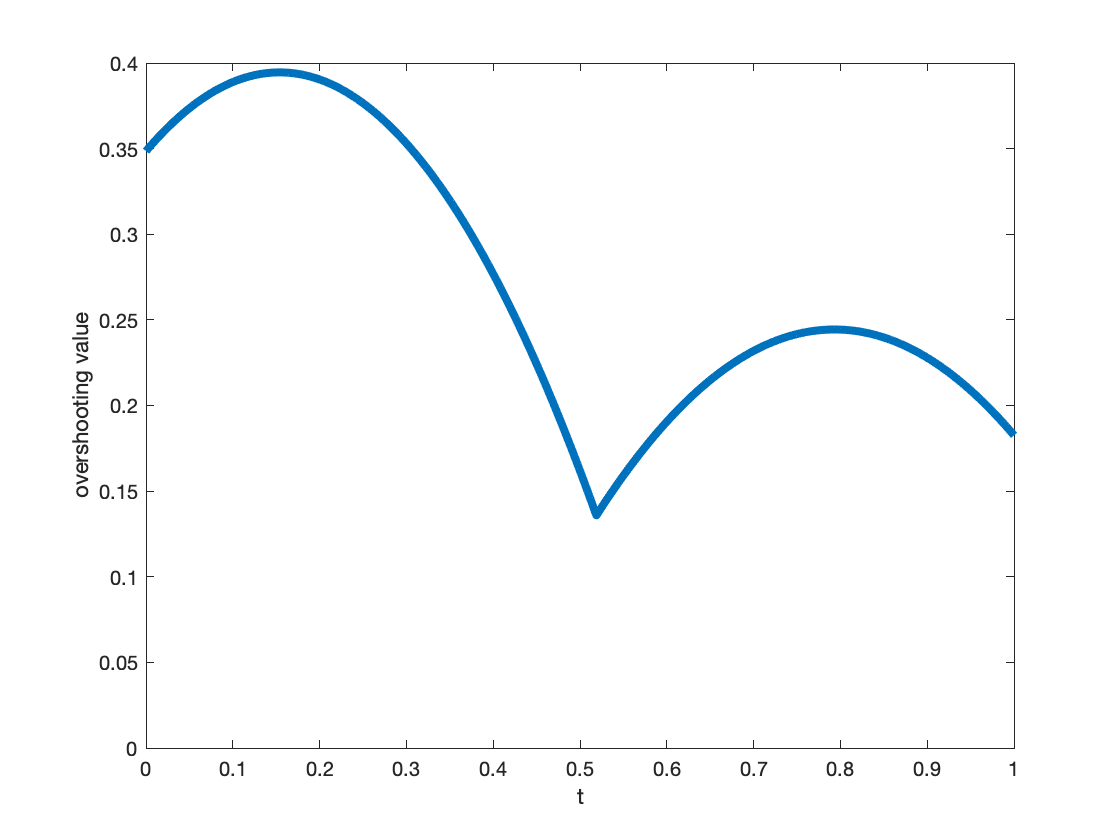}}
\caption{Continuous piecewise linear approximation on a discontinuity non-aligned  nonuniform mesh}
 \label{CP1_non2}
\end{figure}
We have the following $L^2$-error:
\beq
\|u-u_{c,1}\|_{0, I_0} = 2
\sqrt{h (20923 + 13173 t - 221541 t^2 + 450576 t^3 - 250668 t^4)/174243  }.
\eeq
For $0\leq t \leq 1$, $\|u-u_{c,1}\|_{0, I_0}$ is between $0.5349 \sqrt{h}$ and $0.6965 \sqrt{h}$.

\paragraph{\bf An adaptive mesh with coarsening}
For the one-dimensional problem, there is one special case that the overshoot phenomenon that can be eliminated.

Consider the following mesh $\cdots, x_{-2},x_{-1},x_0,x_1,x_2,\cdots$, with $x_0=0$, $x_{\pm 1} = \pm h$, and $x_{\pm 2} = \pm c h$, with $h>0$ and $c>1$. Similar to the uniform grid case, the projection $u_{c,1}$ on $(x_0,x_2)$ can be written as 
$$
u_{c,1} = U_1 \lambda_1 + (1+\epsilon) \lambda_2, \quad x\in(x_0,x_2).
$$
Computing by projection, we get
$$
U_1 = 1 + \dfrac{1-c\epsilon}{2(1+c)}.
$$
This means when the ratio $c>1$ is big enough, the overshoot value can be reduced to a very small value. For a one-dimensional mesh, this is possible by combining the refining and coarsening, i.e., when an a posteriori error indicator on an element $K$ is big, then the element is refined, while for those elements with small indicators, they are combined with the close elements to form bigger elements.  

We test the case with the following example: the domain is $(-1,1)$ and the initial mesh is $[-1;-1/3; 2/3;1]$. We use the exact $L^2$ error as the error estimator, i.e., $\eta_K = \|u-u_{c,1}\|_{0,K}$, and set the refinement/coarsening criteria as "refine those elements whose $\eta_K > 0.6 \max_{\{T\in\cT_h\}}\eta_T$ and coarsen those elements whose  $\eta_K < 0.3 \max_{\{T\in\cT_h\}}\eta_T$". For this very simple problem, the final mesh only contains three elements, a central very small element around $0$, and two big elements on the left and right, receptively, see Fig. \ref{CP1_coarsen}. After 14 iterations, the final mesh is $ [-1; -0.0000814;   0.0000407;   1]$; the numerical solutions at four nodes are $[-1; -1.00008;   0.99998; 1.0]$. The overshoot phenomenon is almost invisible.
\begin{figure}[!htb]
\centering 
~
\subfigure[overshoot values as $t \in (0,1)$]{
\includegraphics[width=0.45\linewidth]{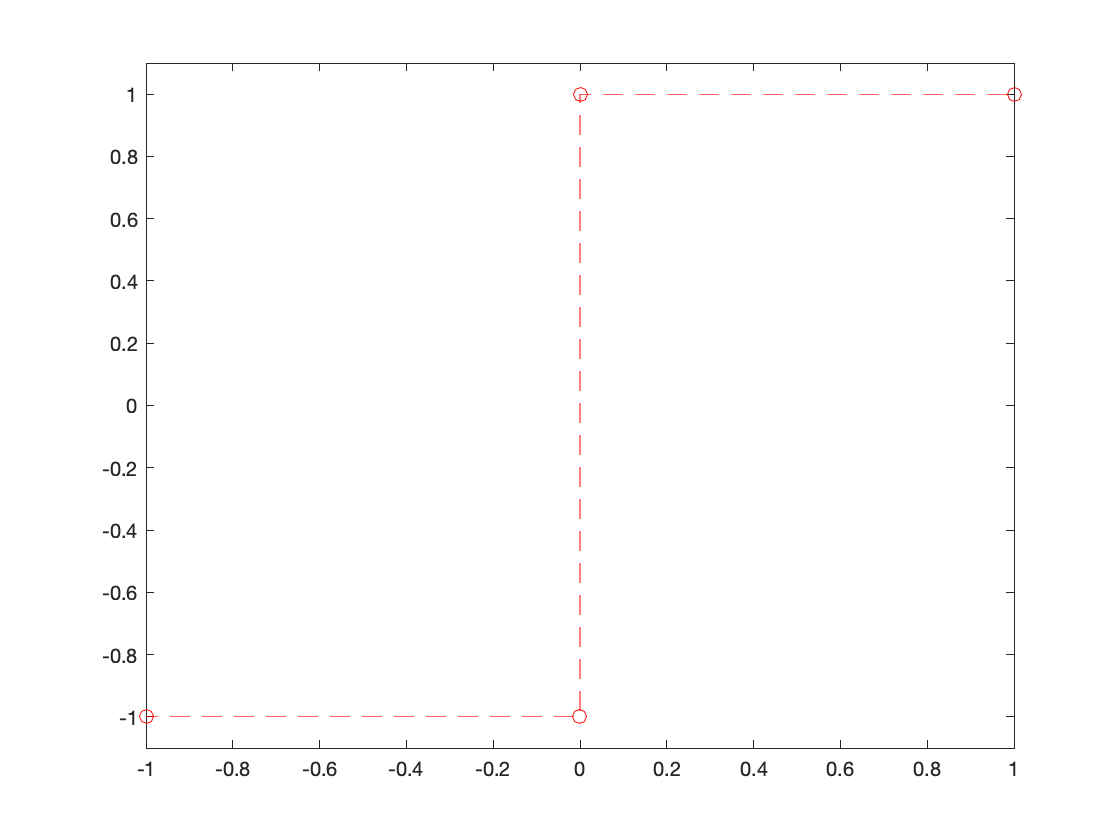}}
\caption{Continuous piecewise linear approximation with adaptive refinement/coarsening}
 \label{CP1_coarsen}
\end{figure}

Note that such case will only happen in a one-dimensional setting. In two and three dimensions, to keep the mesh conforming and regular, the radio of mesh sizes between two adjacent elements sharing a same edge/face cannot be big, thus there will always be non-trivial overshoots. Even for one-dimensional problems, such a high ratio mesh is probably not possible due to other approximation errors to keep the mesh from coarsening, see a numerical test in Section 4.

\section{Numerical Test in 1D: a 1D Singularly Perturbed Problem}
\setcounter{equation}{0}
Consider the following problem:
\beq
 - \epsilon u'' + u = f, \quad \mbox{on  } (0,1), \quad \mbox{and}\quad u(0)=u(1) =0,
\eeq
where 
$$
f = \left\{\begin{array}{cc}
2x & x<1/2, \\
2x-2 & x>1/2.
	\end{array}\right. 
$$
We choose $\epsilon = 10^{-16}$. For $\epsilon$ this small, the solution is almost identical to $f$, and has an extremely sharp layer at $x=1/2$. For the case that mesh size is not small enough to match with the layer, $x=1/2$ can be viewed as a discontinuous location.

We test the problem with three numerical methods: the linear $C^0$-conforming finite element method, linear discontinuous Galerkin finite element method, and the lowest-order mixed finite element method.

Let $0=x_0 <x_1 <\cdots<x_N =1$ be a one-dimensional mesh $\cT_h$ of $(0,1)$, with $K_i =(x_i, x_{i+1})$ and $h_i= x_{i+1}-x_i$. 

\noindent
{\bf Linear conforming finite element method (P1-conforming):} Find $u_h \in S_{1,0}(\cT_h)$, such that
\beq
(\epsilon u_h', v')+(u_h,v) = (f,v), \quad v \in S_{1,0}(\cT_h).
\eeq
where $S_{1,0}(\cT_h) = S_1(\cT_h) \cap H_0^1(\O)$.
Note that for a very small $\epsilon$, this is the $L^2$-projection to $S
_{1,0}(\cT_h)$ we discussed.


\noindent {\bf Linear discontinuous Galerkin finite element problem (P1-DG):} Find $u_{dg} \in P_1(\cT_h)$, such that
\beq
	a_{dg}(u_{dg},v) = (f,v),\quad \forall \, v\in P_1(\cT_h),
\eeq
where 
\begin{eqnarray*}
	a_{dg}(w,v) &:=& (\epsilon w',v')_{\cT_h}+(w,v)_{\cT_h} 
		- \sum_{i=1}^{N-1}\left(\{\epsilon w'(x_i)\}\jump{v(x_i)}+\{\epsilon v'(x_i)\}\jump{w(x_i)}\right) \\
		&& + \sum_{i=1}^{N-1} \dfrac{\epsilon \mu}{h_i+h_{i-1}}\jump{w(x_i)}\jump{v(x_i)},
\end{eqnarray*}
where $\mu>0$ is a big enough number (we can safely choose $\mu=10$).
For a very small $\epsilon$, this can be viewed as the $L^2$-projection to $P_1(\cT_h)$.

To introduce the mixed method, we let $\sigma = -\epsilon u' \in H^1(0,1)$, then $\sigma' +u =f$. The mixed variational formulation is: find $(\sigma, u)\in H^1(0,1)\times L^2(0,1)$, such that,
\beq
\left\{\begin{array}{lclll}
(\epsilon^{-1}\sigma, \tau) - (\tau', u) &=&0,  & \forall \tau \in H^1(0,1),\\[2mm]
-(\sigma', v) - (u,v) & =& -(f,v), &  \forall v\in L^2(0,1).
	\end{array}\right.
\eeq
We choose the approximation spaces to be $S_1$ and $P_0$ for $\sigma$ and $u$, respectively.

\noindent {\bf Mixed finite element problem (P0-mixed):} Find $(\sigma_h, u_h)\in S_1(\cT_h)\times P_0(\cT_h)$, s.t.,
\beq
	\left\{\begin{array}{lclll}
	(\epsilon^{-1}\sigma_h, \tau) - (\tau', u_h) &=&0, &\forall \tau \in S_1(\cT_h), \\[2mm]
-(\sigma_h', v) - (u,v) & =& -(f,v), &  \forall v\in P_0(\cT_h).
	\end{array}\right.
\eeq
For a very small $\epsilon$, this can be viewed as the $L^2$-projection to $P_0(\cT_h)$.

In the computation, we use the robust error estimator for the conforming finite element method in \cite{Ver:98,Kun:02} to drive the adaptive mesh refinement and compute the conforming, DG, and mixed finite element solutions on the same adaptive mesh.

\begin{figure}[!htb]
\centering 
\subfigure[solutions with a discontinuity matched mesh]{ 
\includegraphics[width=0.45\linewidth]{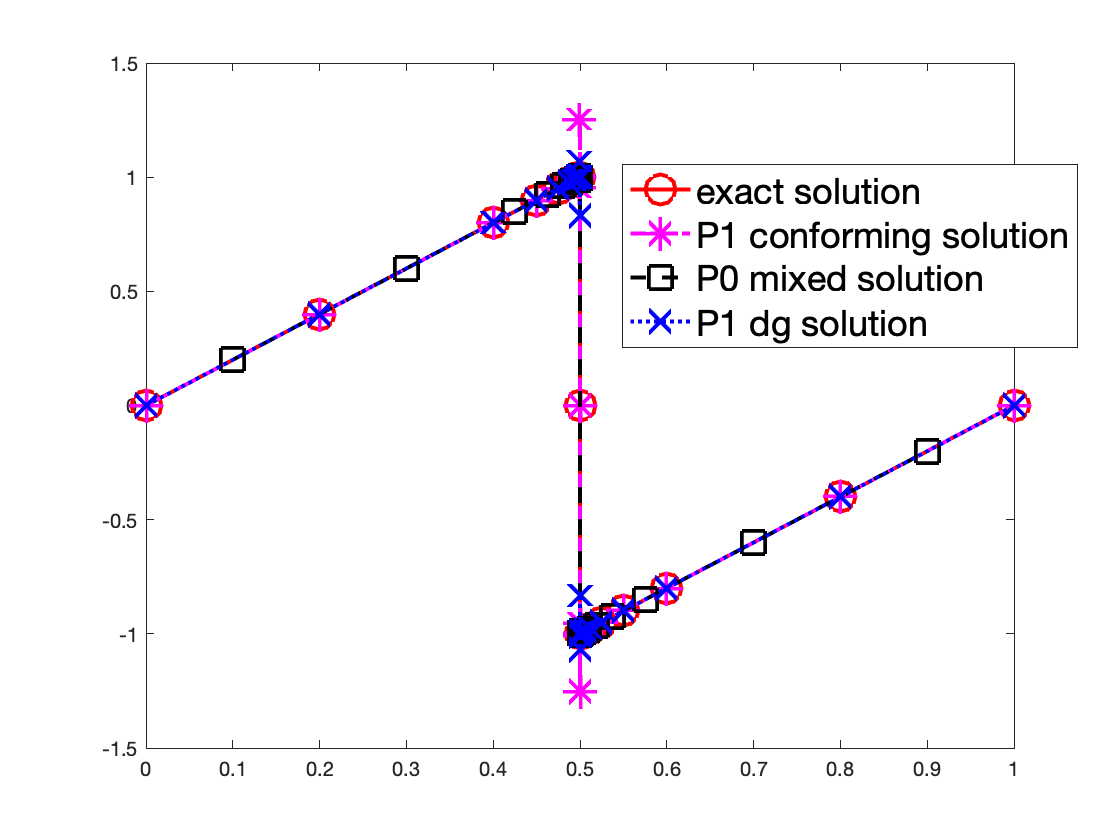}}
~
\subfigure[solutions with a discontinuity non-matched mesh]{
\includegraphics[width=0.45\linewidth]{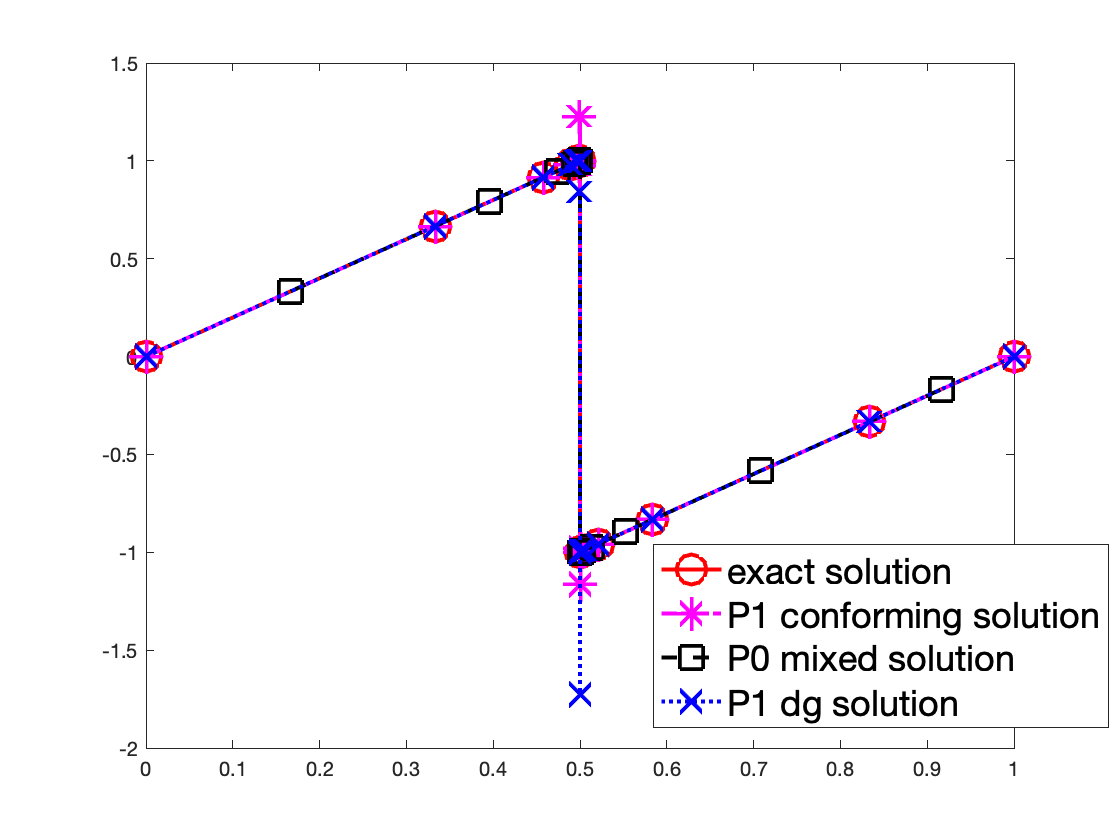}}
\caption{1D reaction-diffusion problem with adaptive refinement on discontinuity aligned and non-aligned meshes}
 \label{mixed_sol}
\end{figure}

\noindent {\bf Discontinuity-aligned mesh case.}
An initial mesh $ [0;0.2;0.4;0.6;0.8;1]$ is chosen. We refine it 20 times with the maximum refinement strategy $\theta = 0.8$. The center nodes are   $0.5+$  $10^{-6}\times [-0.3815;    -0.1907; 0;    0.1907;    0.3815]$. Numerical approximations with linear conforming, P0-mixed, and P1-DG approximations are shown in Fig. \ref{mixed_sol}. The overshoot value for the linear conforming finite element approximation is    $1.2546-1 = 0.2546$, matching the discontinuity match case. The P0-mixed solution has no overshoot because it is a piecewise constant approximation, and the P1-DGFEM has no overshoot since the discontinuity is matched with the mesh.

\noindent {\bf Discontinuity non-aligned mesh case.}
Next, we choose an initial mesh $ [0;1/3;5/6;1]$. With bisections, the discontinuity location will never be exactly matched. The center nodes are   $0.5+$  $10^{-5}\times [
   -0.2543;
   -0.0636;
   -0.0159;$
   $0.0318;
    0.1272;
    0.5086]$.  The values of linear conforming solution $u_h$ at two nodes adjacent to $1/2$ are $1.2284$ and $-1.1647$.  This also matches the discussion in \eqref{os_cp1}. A value of $-1.7242$ is found of P1-DGFEM approximation at the node $x=0.0318$. The value is a little big larger than the discussion in the previous section, where a continuous piecewise linear function is used to approximate a discontinuous step function, while here it is a discontinuous piecewise linear function instead of a step function.  The discontinuous mixed solution has no overshoot in this case.

For a similar problem, on the bottom page 165 of \cite{FR:11}, the authors claimed that the DG method for this problem will have no oscillations. This is imprecise: the DG method will have oscillations unless the discontinuity is matched with the mesh. The numerical test in   \cite{FR:11} uses a uniform mesh with $h=1/N$, which is exactly the discontinuity-aligned case. Thus, the non-overshoot in \cite{FR:11} is due to the use of DGFEM and because of using DGFEM on a discontinuity-aligned mesh.

\noindent {\bf Adaptive mesh with both refinement and coarsening.}
Next, we test the problem with both refinement and coarsening. The initial mesh is chosen to be a uniform mesh with $h =1/17$, and the nodes on the final mesh are $0, 
    0.49262,
    0.49264,
    0.9412,
    1.0000$.
    The $C^0$ conforming finite element solution $u_h$ at these nodes are $
         0,    0.95,   -1.05,   -0.10,         0$ with almost no overshoot, see      Fig. \ref{rd-coarsen}.

\begin{figure}[!htb]
\centering 
\includegraphics[width=0.45\linewidth]{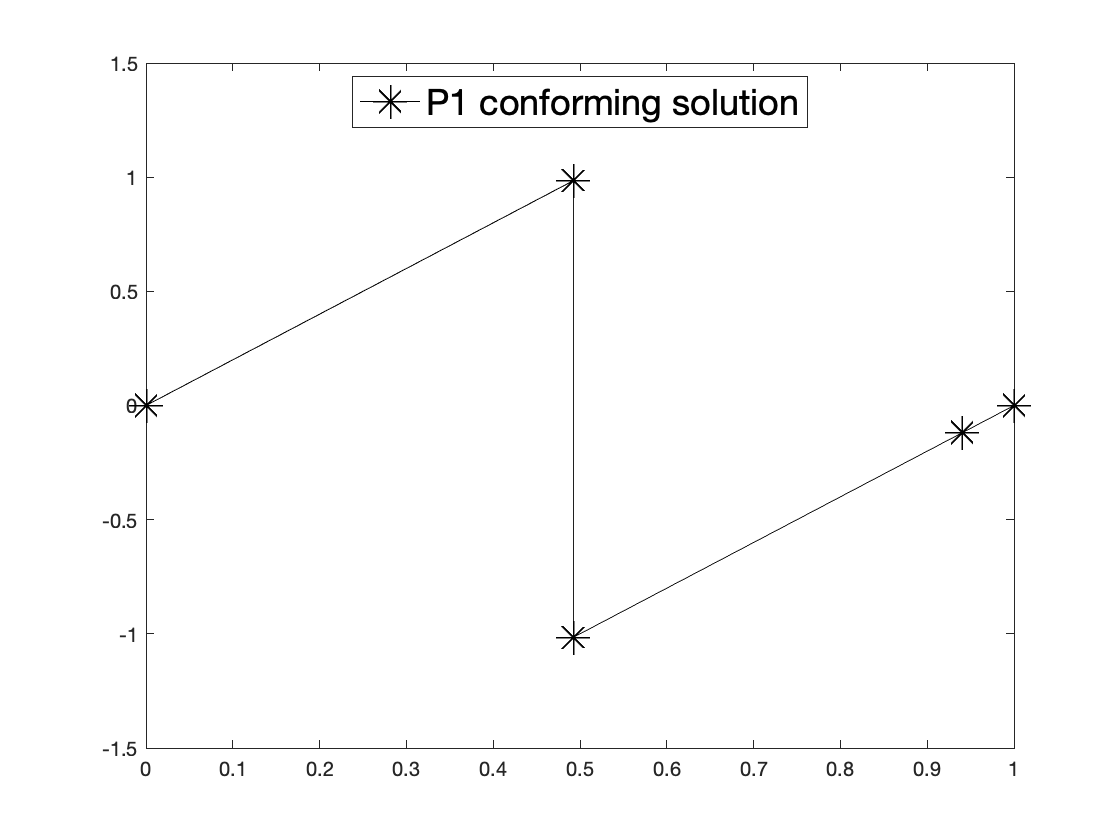}
\caption{1D reaction-diffusion problem with adaptive refinement/coarsening}
 \label{rd-coarsen}
\end{figure}

We should point out that it is not always possible to get a high ratio mesh such as that in Fig. \ref{rd-coarsen} with a  combination of adaptive refinement and coarsening. We test the problem with the right-hand side:
$$
f_2 = \left\{\begin{array}{lc}
x^2 & x<1/2, \\
x^2-1 & x>1/2.
	\end{array}\right. 
$$
Note that there is some approximation error on all elements for the linear finite element approximation. Thus these elements are not coarsened, we eventually have a standard adaptively refined mesh, and an overshoot phenomenon appears, see Fig. \ref{f2}.

\begin{figure}[!htb]
\centering 
\includegraphics[width=0.45\linewidth]{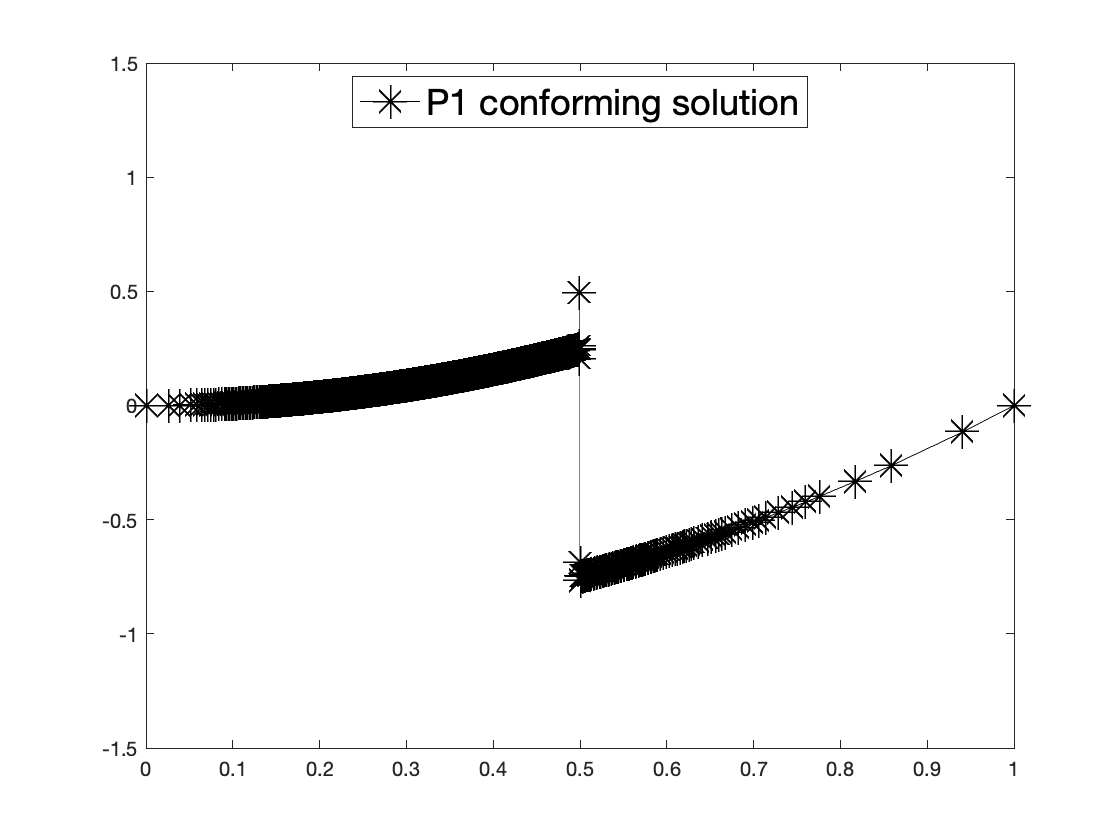}
\caption{1D reaction-diffusion problem, adaptive refinement and coarsening with $f_2$}
 \label{f2}
\end{figure}

\section{Numerical Test in 2D: DG methods  for a Linear Transport Equation}
\setcounter{equation}{0}
Consider the following linear transport (advection) equation:  
\begin{align} \label{transporteqn}
\gradt(\bbeta u)+\gamma u &= f \quad \mbox{in} \,\ \O, \\ \nonumber
                        u &= g \quad \mbox{on} \,\ \Gamma_{-},
\end{align}
where the advective velocity field $\bbeta=(\beta_1,\beta_2)^{T}\in [C^1(\overline{\O})]^2 $ is a
vector-valued function defined on $\bar{\O}$ and $\gamma \in L^{\infty}(\O)$. The inflow part of $\p \O$, $\Gamma_{-}=\{x\in\p \O:\bbeta(x)\cdot\bn(x)<0\}$,  is defined in the usual fashion, where $\bn(x)$ denotes the unit outward normal vector to $\p\O$ at $x\in\p\O$. 

The DGFEM for the transport equation is \cite{BMS:04}: find $u_h \in P_k$, $k=0$ or $1$ such that
\beq
a_{dg}(u_h,v_h) = f(v_h), \quad \forall v_h \in P_k,
\eeq
where the bilinear form and linear form are defined by 
\begin{eqnarray*}
a_{dg}(w_h,v_h) &:=& \sum_{K\in \cT_h} (w_h, -\bbeta \cdot \nabla v_h+\gamma v_h)_K + \sum_{F \in \cE_I \cup \cE_{out}} (\{\bbeta\cdot\bn w_h\}_{up}, \jump{v_h})_F, \\
f(v_h) &=& (f, v_h)- \sum_{F\in \cE_{in}} (\bbeta\cdot\bn g, v_h)_F,
\end{eqnarray*}
where the term $\{\bbeta\cdot\bn w_h\}_{up}$ is the usual upwind flux, meaning that $\{ \bbeta\cdot\bn w_h\}_{up}$ takes the value on the inflow side of the edge/face $F$. The sets $\cE_{in}$, $\cE_{I}$, and $\cE_{out}$ are the sets of inflow, interior, and outflow edges, respectively. We call the methods P0-DGFEM and P1-DGFEM for $k=0$ and $1$.
We use a residual type of error estimator similar to that developed in \cite{Burman:09} to drive the adaptive mesh refinement.

\subsection{A discontinuous solution on a non-matching mesh}
Consider the problem: $\O = (0,2)\times (0,1)$ with $\bbeta = (0,1)^T$. 
The inflow boundary is $\{x\in (0,1), y=0\}$. Let $\gamma =0$ and $f=0$. Choose the inflow boundary condition such that the exact solution is
$$
u(x,y) = \left\{ \begin{array}{lll}
0 & \mbox{if} & x< \pi/3, \\[2mm]
1 & \mbox{if} & x > \pi/3.
\end{array} \right.
$$
The initial mesh is shown on the left of Fig. \ref{inimesh_pi}. The point $(\pi/3,0)$ and $(1,1)$ are the bottom central node and the top central node, respectively. Although the inflow boundary condition is matched, the mesh and its subsequent bisection mesh are not aligned with the discontinuity. 

In the center of Fig. \ref{inimesh_pi}, a final adaptive mesh is shown with many refinements along the discontinuity since the mesh is not matched. On the right of Fig. \ref{inimesh_pi}, we show the overshoot values of the P1-DGFEM on different levels of refined meshes, with the overshoot value defined as:
\beq \label{ooss}
{\tt os} = \max_{\bx\in \cN}\{u_h(\bx)-1, -u_h(\bx), 0\},
\quad
\cN \mbox{  is the collection of nodes of the finite element mesh}.
\eeq

From the figure, we see that the overshoot values oscillate between 0.15 and 0.35 and never go to zero when the mesh is adaptively refined, which is matched with our analysis. For the P0-DGFEM, the overshoot value is in the order of machine accuracy, $10^{-16}$, so there is no overshoot at all.

\begin{figure}[!htb]
\centering 
\subfigure[initial mesh]{ 
\includegraphics[width=0.3\linewidth]{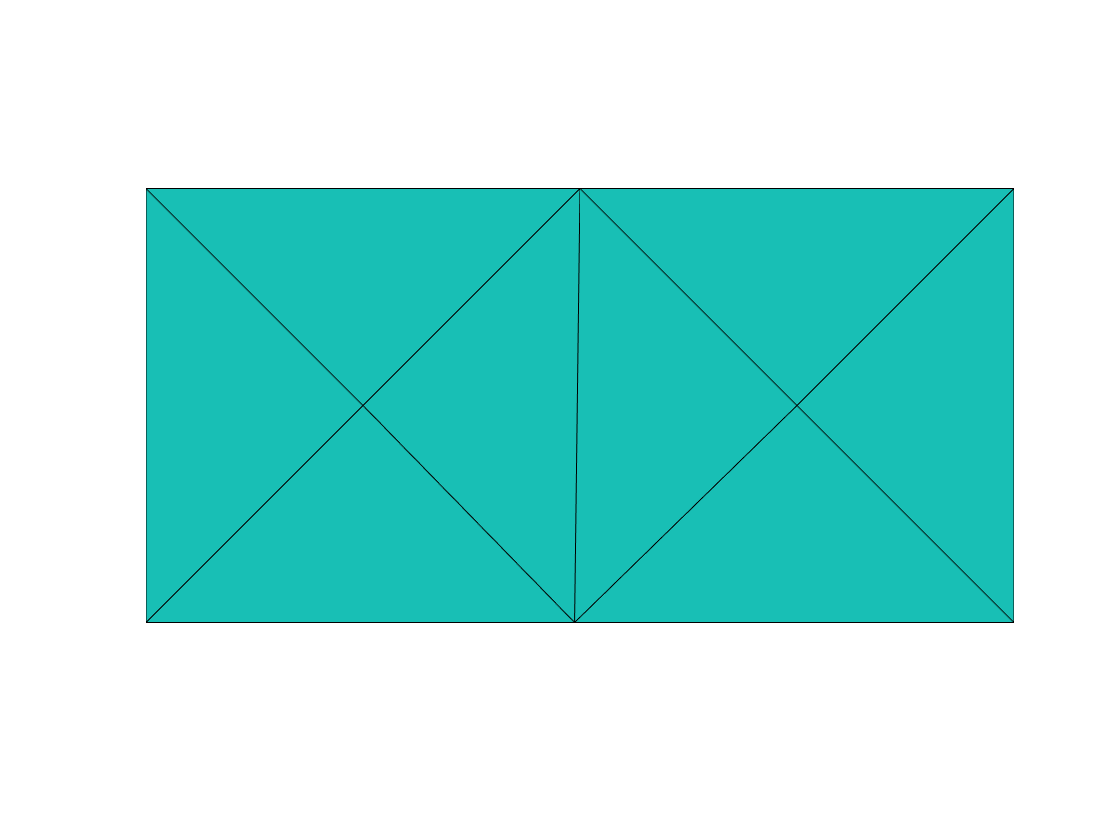}}
~
\subfigure[final adaptive mesh]{ 
\includegraphics[width=0.3\linewidth]{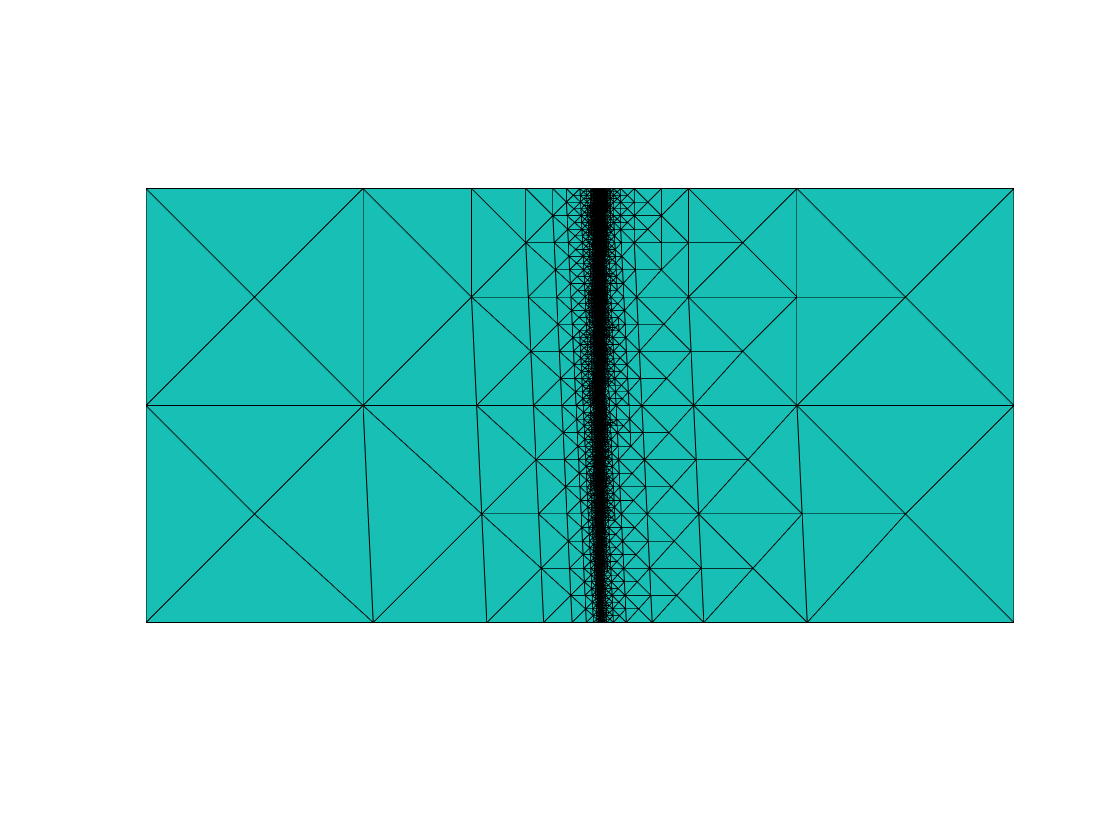}}
~
\subfigure[overshoots on different levels of refined meshes with P1-DGFEM ]{
\includegraphics[width=0.3\linewidth]{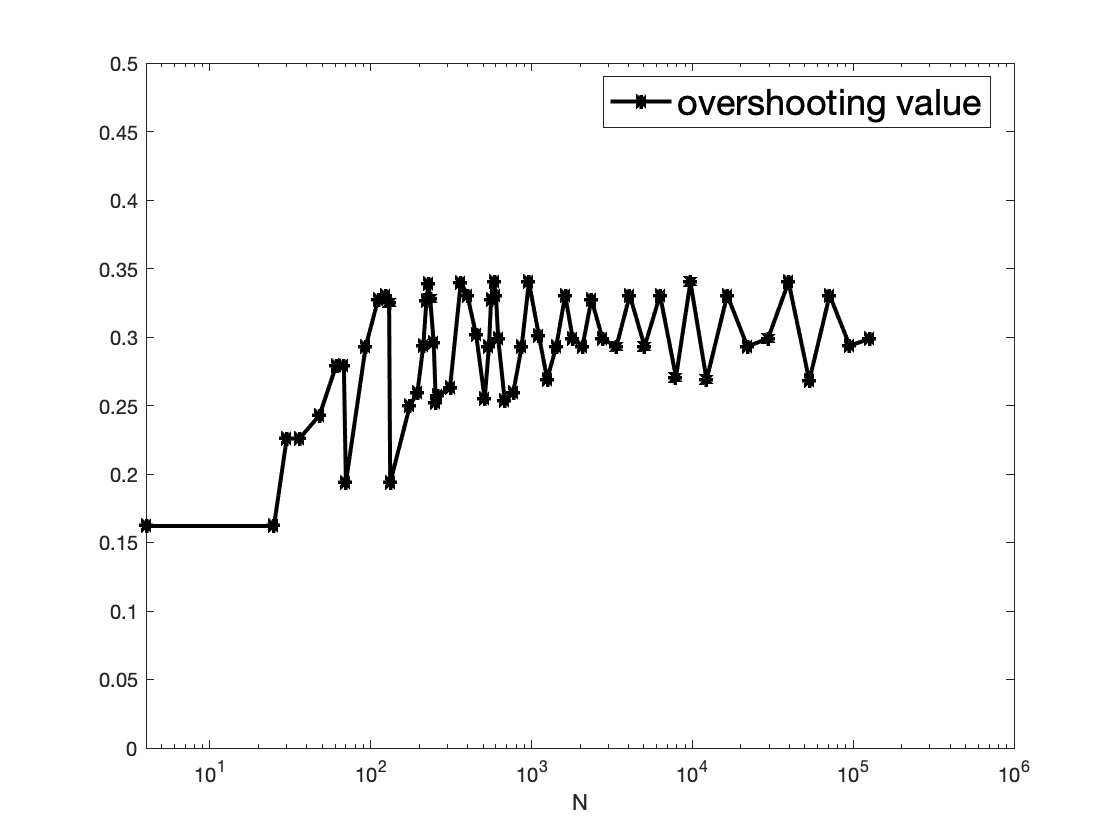}}
\caption{Piecewise constant solution transport problem with a non-matching grid test problem}
 \label{inimesh_pi}
\end{figure}

Since the exact solution of the problem is not changing with respect to the $y$-coordinate, we plot a projected solution by plotting the numerical solutions of the element center ($P_0$)/ nodes ($P_1$) suppressing the $y$-coordinate.  In Fig. \ref{projection_PWCNM}, we show the P0-DGFEM projected solution computed on the final adaptive mesh on the left and the P1-DGFEM one on the right. We clearly see the overshoot of P1-DGFEM.  

\begin{figure}[!htb]
\centering 
\subfigure[P0-DGFEM]{ 
\includegraphics[width=0.45\linewidth]{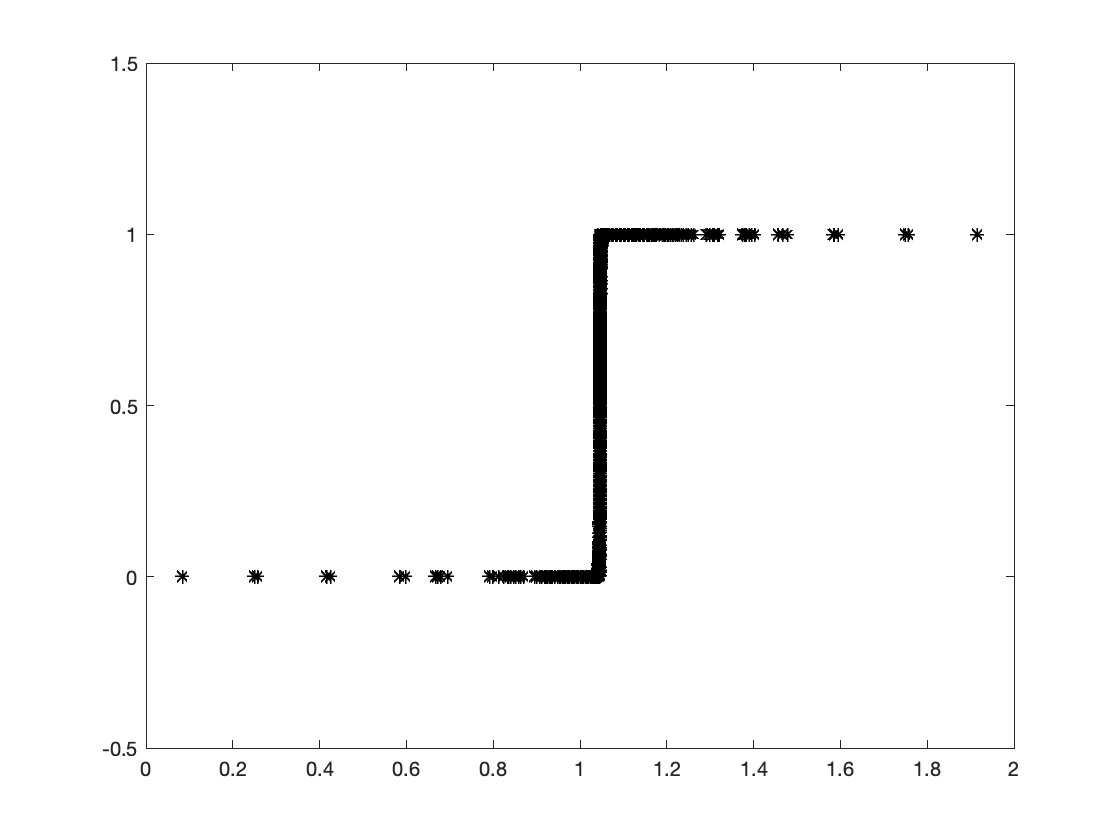}}
~
\subfigure[P1-DGFEM]{
\includegraphics[width=0.45\linewidth]{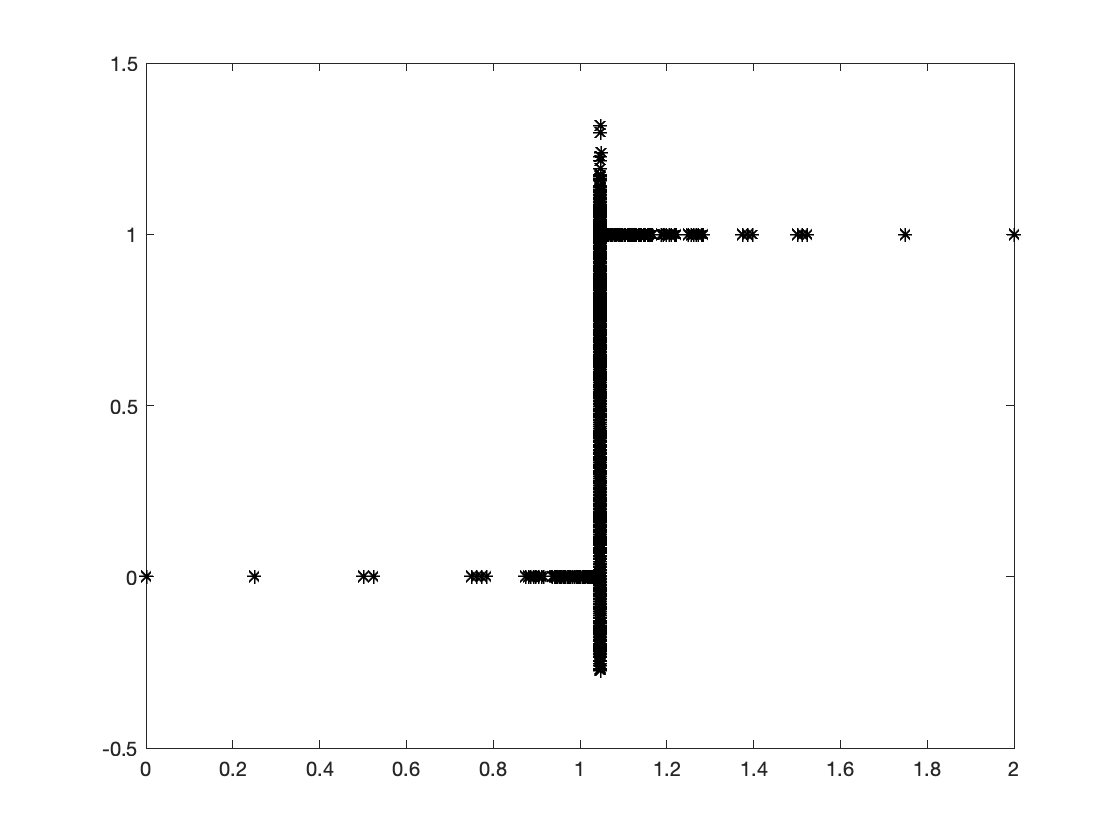}}
\caption{Projected solutions with P0- and P1-DGFEMs for piecewise constant solution transport problem on a non-matching mesh}
 \label{projection_PWCNM}
\end{figure}


\subsection{Curved transport problem 1}
Consider the problem on the half disk  $\O = \{(x,y) \colon x^2+y^2<1; y>0\}$. The inflow boundary is $\{-1<x<0; y=0\}$. Let $\bbeta = (\sin \theta, -\cos \theta)^T = (y/\sqrt{x^2+y^2}, -x/\sqrt{x^2+y^2})^T$, with $\theta$ is the polar angle. Let $\gamma=0$, $f=0$, and the inflow condition and the exact solution be
$$
g = \left\{ \begin{array}{lll}
1 & \mbox{if}  &-1<x<-0.5, \\[2mm]
0 & \mbox{if}  &-0.5<x<0,
\end{array} \right.
\mbox{and}\quad
u = \left\{ \begin{array}{lll}
1 & \mbox{if  }  x^2+y^2 > 0.25, \\[2mm]
0 & \mbox{otherwise}. 
\end{array} \right.
$$ 

\begin{figure}[!htb]
\centering 
\subfigure[initial mesh]{ 
\includegraphics[width=0.3\linewidth]{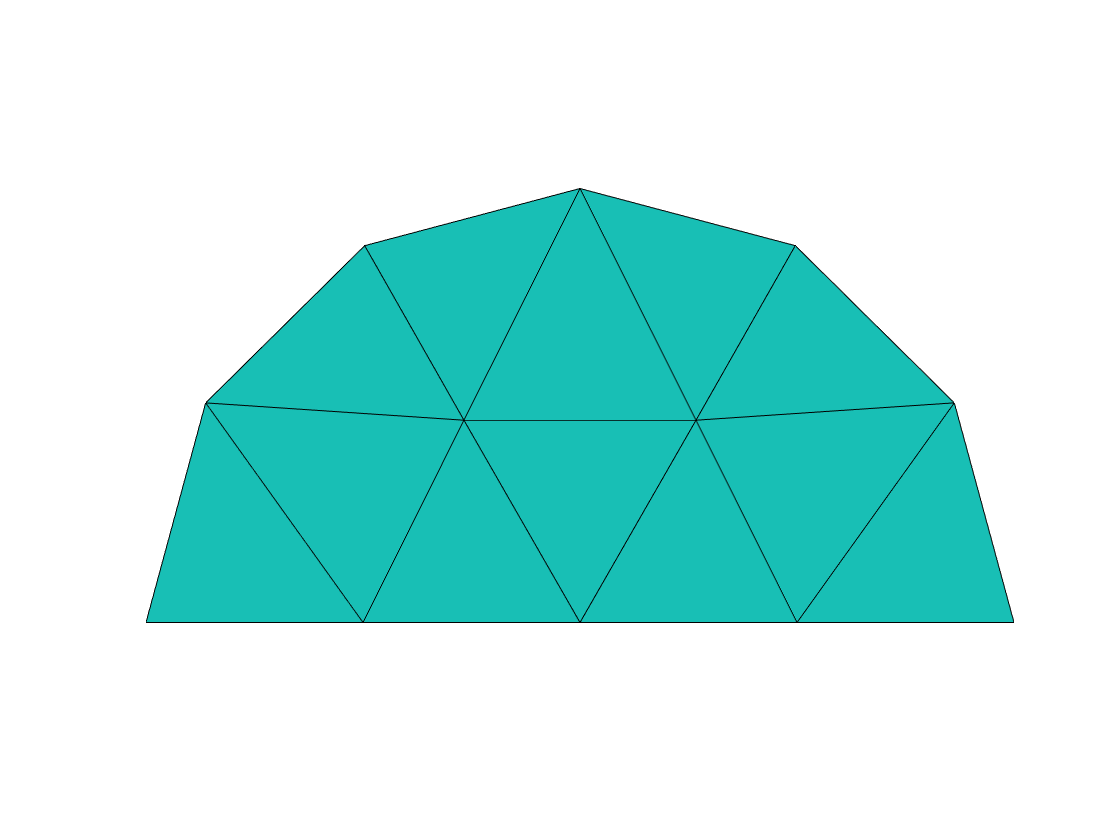}}
~
\subfigure[adaptive mesh]{ 
\includegraphics[width=0.3\linewidth]{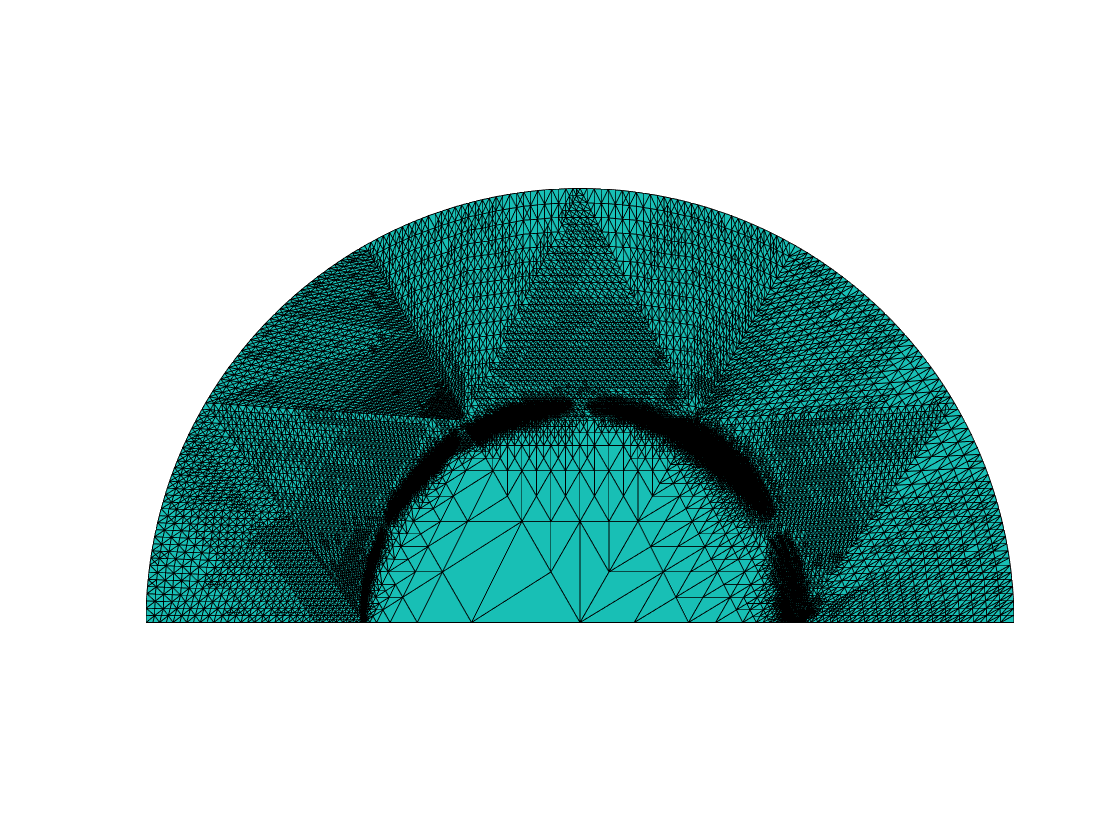}}
~
\subfigure[os by P1-DGFEM]{
\includegraphics[width=0.3\linewidth]{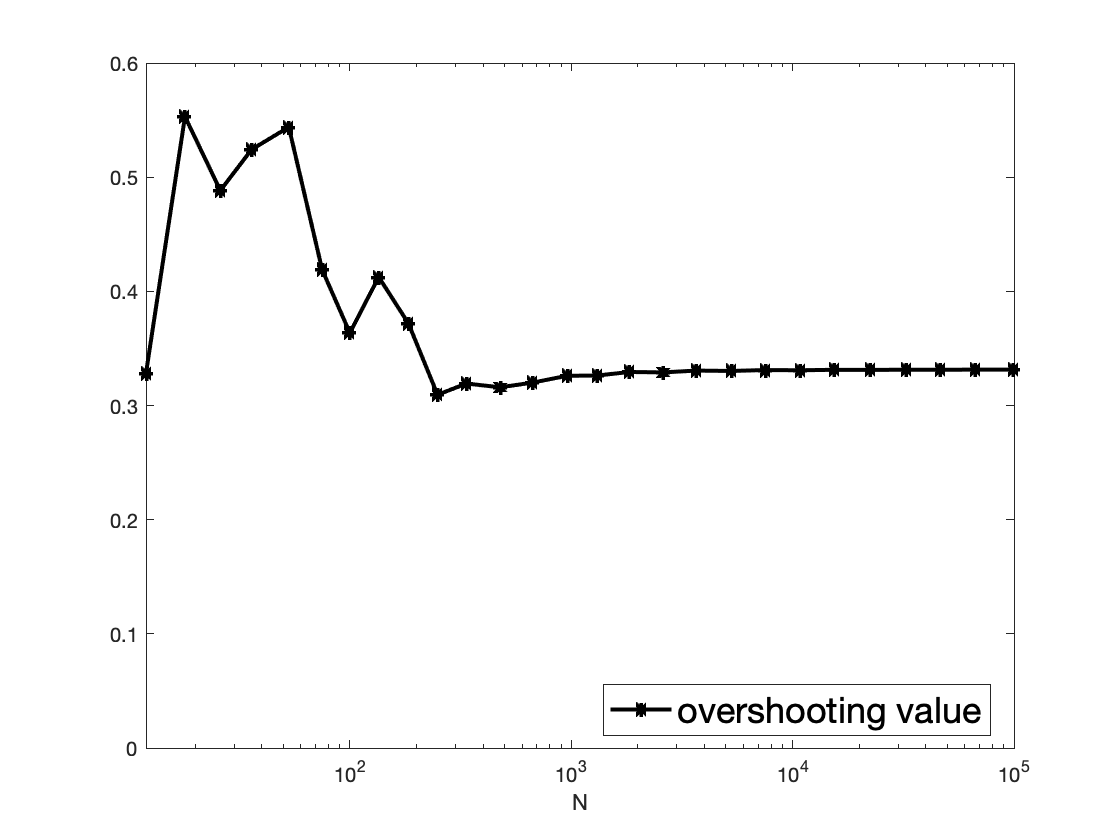}}
\caption{P1-DGFEMs for the curved transport problem 1}
 \label{os_curved}
\end{figure}

\begin{figure}[!htb]
\centering 
\subfigure[P0-DGFEM]{ 
\includegraphics[width=0.45\linewidth]{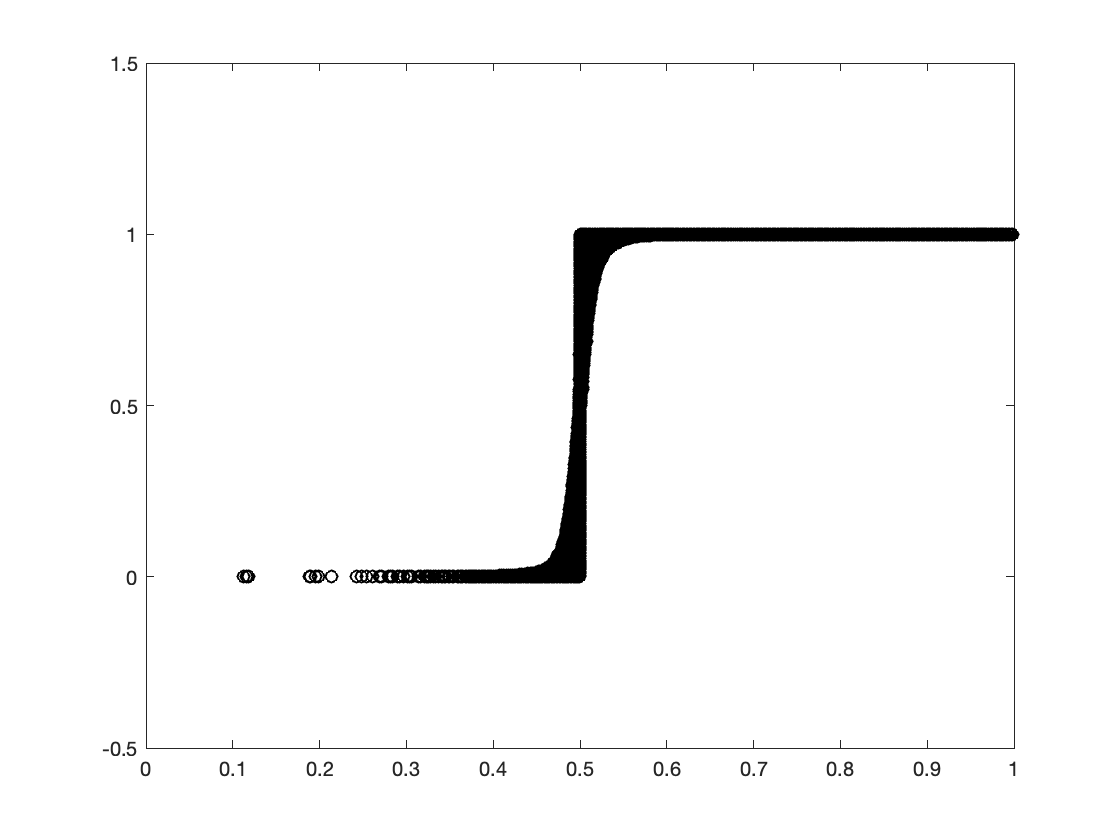}}
~
\subfigure[P1-DGFEM]{
\includegraphics[width=0.45\linewidth]{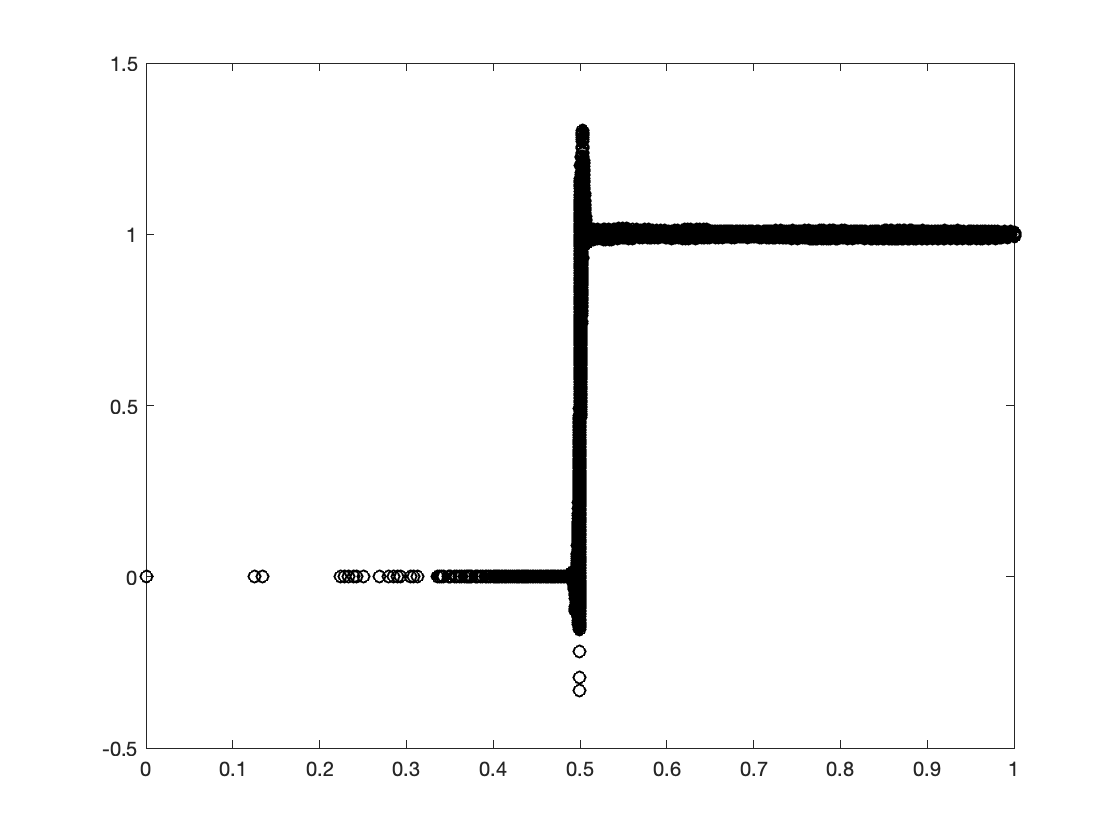}}
\caption{Projected solutions with P0 and P1-DGFEMs for curved transport problem 1}
 \label{projection_curved}
\end{figure}

The initial mesh  is shown on the left of Fig. \ref{os_curved}. On the center of Fig. \ref{os_curved}, a final adaptive mesh is shown. On the right of Fig. \ref{inimesh_pi}, we show the overshoot values of the P1-DGFEM, with the overshoot value \eqref{ooss}.

From the figure, we see that the overshoot values oscillate between 0.3 and 0.6 and never go to zero with the mesh refinements.

We plot projected solutions by plotting the numerical solutions of the element center ($P_0$)/ nodes ($P_1$) with respect to the radius in Fig. \ref{projection_curved}, we show the P0-DGFEM projected solution computed on the final adaptive mesh on the left and the P1-DGFEM one on the right. We clearly see the overshoot of the P1-DG approximation but the non-overshoot of the P0-DG approximation.

\subsection{Curved transport problem 2}
Consider the following problem: $\O = (0,1)^2$ with $\bbeta = (y+1,-x)^T/\sqrt{x^2+(y+1)^2}$, $\gamma=0.1$, and $f=0$. The inflow boundary is $\{x=1, y\in (0,1)\} \cup \{x\in (0,1), y=0\}$, i.e., the west and north boundaries of the domain.  Choose $g$ such that the exact solution $u$ is
$$
u = \dfrac{1}{4}\exp\left(\gamma r\arcsin \left(\dfrac{y+1}{r}\right)\right) 
\arctan \left(\dfrac{r-1.5}{\epsilon}\right), \quad \mbox{with}\quad r = \sqrt{x^2+(y+1)^2}.
$$
When $\epsilon = 10^{-10}$, the layer is never fully resolved in our experiments and can be viewed as discontinuous. In Fig. \ref{burman}, we show numerical solutions by P0- and P1-DGFEMs. It is easy to see that there is no overshoot for the P0-DGFEM solution, and a non-trivial overshoot can be found for the P1-DGFEM solution. 

\begin{figure}[!htb]
\centering 
\subfigure[P0-DGFEM]{ 
\includegraphics[width=0.45\linewidth]{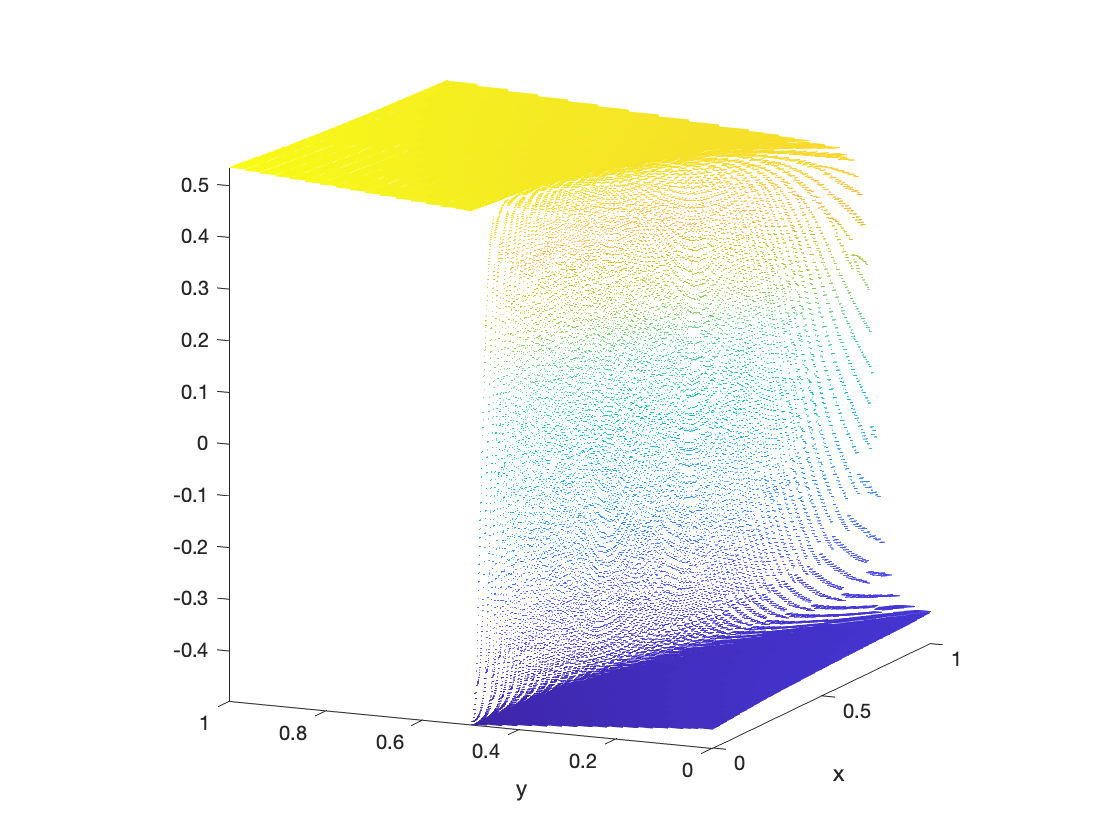}}
~
\subfigure[P1-DGFEM]{
\includegraphics[width=0.45\linewidth]{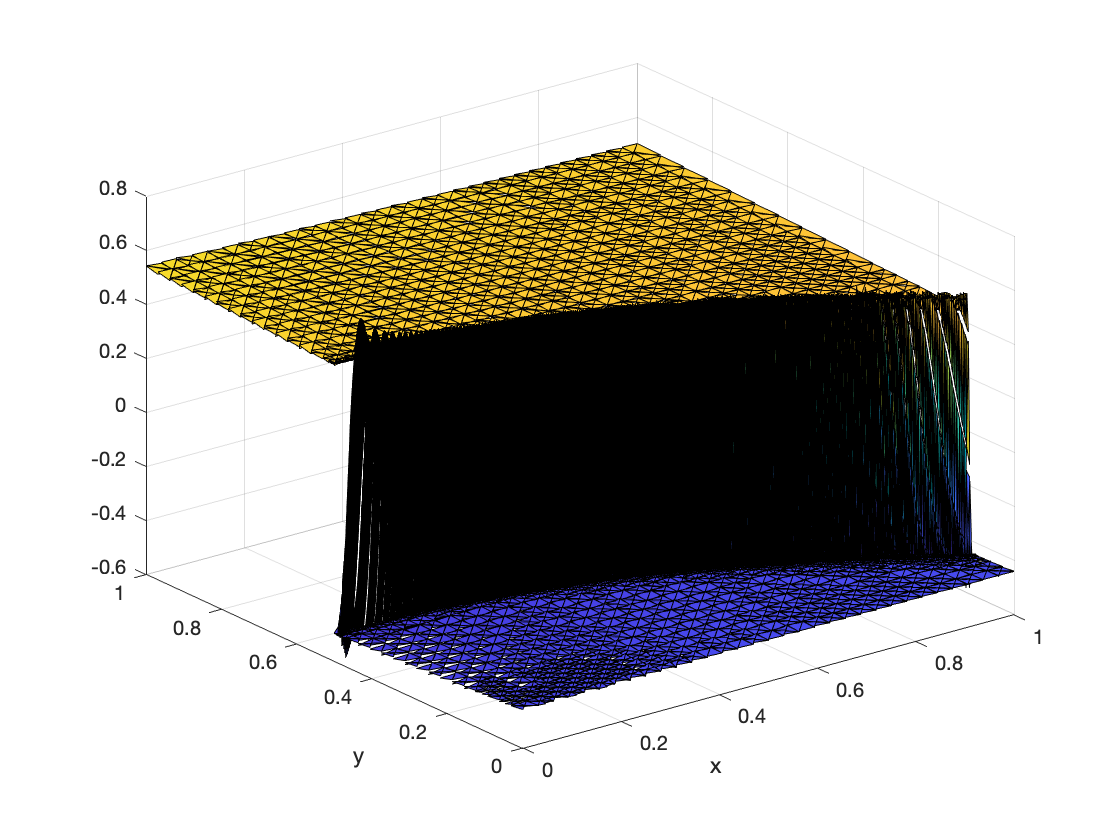}}
\caption{Numerical solutions with adaptive P0- and P1-DGFEMs for Curved Transport Problem 2}
 \label{burman}
\end{figure}

Similar results are available for these three numerical tests with $H(\divvr)$ flux-based least-squares finite element formulations can be found in \cite{LZ:18,LZ:19}.

\section{Some Discussions of Numerical Tests in the Literature}
In \cite{ACFFJLU:11},  various numerical methods are tested for the convection-dominated convection-diffusion equations Hemker problem, which contains boundary and interior layers. The interior layers are not matched with the mesh. The authors also conclude that the classic finite volume method has minimal overshoots. The classic finite volume method uses piecewise constant approximation. Combined with adaptive methods, it has a hope of getting an accurate and spurious oscillations-free solution.

In \cite{LWY:17}, a least-squares-based weak Galerkin finite element method is applied to reaction-diffusion problems. The discontinuous finite element space is used to approximate the primal variable in the method. In a one-dimensional example (Example 7) with an almost discontinuity at the center $x=0.5$. From Figures 7.6-7.8 of the paper, we can see a grid point on $x=0.5$. This is the discontinuous approximation in a discontinuity matched case, we should not observe the Gibbs phenomenon, and it matches the numerical test in  Figures 7.6-7.8 of \cite{LWY:17}.

In \cite{LYZZ:18}, a weak Galerkin method is introduced to solve convection-dominated convection-diffusion-reaction equations. Its numerical example 4 is $-\epsilon \Delta u + \gradt (\bbeta u) + \gamma u = 0$ with discontinuous boundary conditions. For a very small $\epsilon$ and a mesh of $1024$ uniform squares, it is the linear advection equation \eqref{transporteqn}. In Figs. 7-8, the Gibbs phenomenon is observed near the discontinuity. The authors "{\em believe that a better numerical approximation can be obtained if the Dirichlet boundary conditions are enforced weakly}", while it is most likely due to that the mesh is not aligned with the discontinuity as opposed to \cite{LWY:17}.

In \cite{CK:20}, for the singularly perturbed reaction-diffusion problem, a two-step finite element method is developed. First, the flux is computed independently, then the numerical solution $u_h$ is recovered in $P_k$. In Theorem 6.5 of \cite{CK:20}, $\|u -u_h\|_{L^\infty(K)} \leq  \|u -P_h u_h\|_{L^\infty(K)} + \mbox{extra terms}$ are proved, where $P_h$  is the $L^2$-projection. The numerical results show that when $k=0$ is chosen, then $u_h$ does not show a significant numerical oscillation. This matches our discussion. On the other hand, if $k\geq 1$ is chosen and the mesh is not aligned with discontinuity, then the term $\|u -P_h u_h\|_{L^\infty(K)}$ will be big, and the Gibbs phenomenon will appear.

\section{A Short Discussion of $L^1$-based methods}
\subsection{Discontinuous Approximations}
For the discontinuous function $u$ defined in \eqref{u} and an interval $I$, consider the following $L^1$ minimization problem
\beq
\mbox{Find}\quad w_h\in V_h,\quad \mbox{such that } \|u-w_h\|_{L^1(I)} \rightarrow \min !
\eeq
\subsubsection{Discontinuous piecewise constant approximations}
In this simple case, we only need to compute a constant $c$ on the interval $I_0$. We should discuss two cases, $|c|> 1$ and $|c|\leq 1$. It is easy to see that the minimizer should happen when $|c|\leq 1$. 

Then 
\beq
\|u-c\|_{L^1(I_0)} = \int_{-th}^{(1-t)h} |u-c| dx = h(1-c + 2tc).
\eeq
After a simple discussion of $t$ in two cases, we have:
\begin{eqnarray*}
&&\mbox{if } 1/2\leq t\leq 1, \quad \mbox{then } c=-1 \mbox{ and } \|u-c\|_{L^1(I_0)}= 2(1-t)h,\\
&&\mbox{if } 0\leq t\leq 1/2, \quad \mbox{then } c=1 \mbox{ and } \|u-c\|_{L^1(I_0)}= 2th.
\end{eqnarray*}
Since $c$ is $1$ or $-1$, then there is no overshoot.

\subsubsection{Discontinuous Piecewise Linear Approximations}
Since the $L^1$-minimization problem is more complicated than a simple $L^2$-projection, we only discuss a simple case to show that $L^1$-minimization can cause an overshoot with a $P_1$ approximation.

For the interval $I_0= [-th, (1-t)h]$, we consider the simplest and symmetric case $t=1/2$. As discussed in Section \ref{S322}, let $w_h = W_{-1}\lambda_{-1} + W_{1}\lambda_{1}$. Due to symmetry, we can safely assume that $w_h(0)=0$ and $W_{1} = -W_{-1} \geq 0$; thus, we only need to consider the following $L^1$-minimization problem on the positive half of $I_0$ (with $\lambda_1 = 2x/h$):
\beq
\mbox{Find } W_1\geq 0:
\|u-w_h\|_{L^1(I_0)} = 2\int_{0}^{h/2} |u- W_{1}\lambda_1| dx
=h\int_{0}^{1} |1- W_{1}x| dx
 \rightarrow \min!
\eeq
The problem can be explicitly solved by discussing two cases: $W_1\geq 1$ or $0\leq W_1 \leq 1$:
\beq
J(W_1):=\int_{0}^{1} |1- W_{1}\lambda_1| dx = 
\left\{ \begin{array}{lll}
1/W_1+W_1/2 -1 & \mbox{if  }  W_1 \geq 1, \\[2mm]
1-W_1/2  & \mbox{if  }  0\leq W_1 \leq 1, . 
\end{array} \right.
\eeq
A simple calculation shows that 
\beq
\min_{0\leq W_1}J(W_1) = \sqrt{2}-1 \quad \mbox{when}\quad W_1 = \sqrt{2}.
\eeq
This shows a non-trivial overshoot appears for $t=1/2$. For other cases $t\in [0,1]$, the Gibbs phenomenon still can happen (at least for $t$ close to $1/2$). We skip the detailed computation since we usually cannot control discontinuity location for non-matched cases.

\subsection{Continuous Approximations}
\subsubsection{Continuous Piecewise Linear Approximations in 1D}
The discussion of $L^1$-best approximation on 1D can be found in Theorem 1.1 (a boundary discontinuity case) and Theorem 1.4 (a jump discontinuity case) of \cite{HRZ:19b}. The theorem concludes that, under special cases, for example, a uniform grid (or a less restrictive case that the elements cannot be too small compared to their neighboring element closer to the discontinuity) with a {\bf discontinuity matched mesh}, the Gibbs phenomenon can be eliminated. In other cases, the Gibbs phenomenon will happen.   

\subsubsection{Continuous Piecewise Linear Approximations in 2D}

The discussion of $L^1$-best approximation on 2D can be found in Theorem 1.5 (boundary discontinuity case) and numerical tests of \cite{HRZ:19b}. It is shown in \cite{HRZ:19b} that even for the boundary discontinuity case, which is a discontinuity matched mesh, the $L^1$-best approximation can only eliminate the Gibbs phenomenon on certain very special meshes. 

For the general case of discontinuity non-aligned meshes, the situation is more complicated.  From observations of the 1D cases and 2D discontinuity aligned mesh cases in \cite{HRZ:19b}, the Gibbs phenomenon still appears for some $L^1$-best approximations. It is also observed in some numerical tests of \cite{HRZ:19b}, that the mesh adaptivity might be helpful in reducing the overshoot, but the underlying mechanism is unclear.

\section{A Table of Results}

\begin{table}[htp]
\caption{Gibbs Phenomenon in Different Situations (O: overshoot, NO: no overshoot)}
\begin{center}
\begin{tabular}{|c|c|c|c|c|c|c|}
\hline
 mesh 	& \multicolumn{2}{c}{$L^2$-projection} & & \multicolumn{2}{c}{$L^1$-best approximation} &  \\
\hline
			   	& $P_0$ & $P_k$ ($k\geq 1$) & $S_k$ ($k\geq 1$) 	& $P_0$ & $P_k$ ($k\geq 1$) & $S_k$ ($k\geq 1$) \\
\hline
Discontinuity-aligned (include boundary layer) 	& NO & NO & O & NO & NO& both \\
\hline
Discontinuity-nonaligned & NO & O &O& NO & O& O \\
\hline
\end{tabular}
\end{center}
\label{tab}
\end{table}%
In Table \ref{tab}, we list the Gibbs phenomenon in different situations. The computation of an $ L^2$-projection method is relatively easy, a linear equation needs to be solved in an $L^2$ projection. On the other hand, methods based on $L^1$-best approximation methods are more complicated, and iterative optimization algorithms are needed; see for example some algorithms suggested in \cite{Guermond:04}.
 From the comparisons of the table, we do not recommend such $L^1$-based methods due to their similar behavior to $L^2$-based methods and more complicated implementations. 

For $L^2$-minimization-based methods, we have the following conclusions.
For the simple discontinuity-aligned mesh case, piecewise discontinuous approximations are always good. For the general non-matched case, the piecewise discontinuous constant approximation combined with adaptive mesh refinements is a good choice to achieve accuracy without overshoots. For discontinuous piecewise linear approximations, non-trivial overshoots will be observed unless the mesh is matched with discontinuity. For continuous piecewise linear approximations, non-trivial overshoots will always be observed under regular meshes.

\section{Concluding Remarks}
In this paper, we study the behavior of using adaptive continuous and discontinuous finite elements to approximate discontinuous and nearly discontinuous PDE solutions from an approximation point of view.

For a singularly perturbed problem with a small diffusion or other equations with a transient layer, the numerical solution can accurately approximate the PDE solution when the mesh is fine enough around the transient layer.  This is the case that the belief that the adaptive finite element works come from.

For a pre-asymptotic mesh of such layer problems, the solution can be viewed as discontinuous, and we can treat it as a genuinely discontinuous case.

To approximate a discontinuous solution, the global continuous function is always a bad choice since it will always have non-trivial oscillations. It is also hard to impose a discontinuous boundary condition in a continuous approximation space. If the exact location of the discontinuity is known, then a high-order discontinuous approximation on a discontinuity-matched mesh is a good choice. For singularly perturbed problems with small diffusion, this mesh-matched case is reduced to that the mesh interface lies within the small transient layer. If the location of the discontinuity is unknown, the only way to avoid the overshoot is to use piecewise constant approximations for those elements where a discontinuity cuts through. Discontinuous piecewise linear function approximations will have a non-trivial overshoot for most cases. The adaptive finite element method with a good a posteriori error estimator can reduce the $L^2$-norms (or $H^1$ and other integration-based norms) of the error.  However, it cannot reduce the non-trivial oscillations unless a piecewise constant approximation is used in the discontinuity-crossing elements. 

The non-trivial overshoot we find in the paper is of the same magnitude as the Gibbs phenomenon in the Fourier analysis, where a number around $0.17$ of the jump gap is found, see \cite{HH:79}.   

With the above results in mind, an hp adaptivity approach is preferred when designing finite element methods for problems with unknown location discontinuous solutions.  We can use high-order finite elements (continuous or discontinuous) in the smooth region and use piecewise constant approximations in the discontinuity-crossing elements. It is also worth mentioning that the method proposed by Lax \cite{Lax:06} is precisely such a method in the context of the finite difference method.

It is also interesting to compare the result to various {\em order barrier theorems} in the finite difference/finite volume schemes, see \cite{Wesseling:01}. 
These theorems are developed in a different context, with the finite difference like method, and are proved by methods very different from the discussion in this paper. On the other hand, the discussion in this paper is from the approximation point of view without a specific partial differential equation. It is ongoing research to find the connections between these order barrier theorems and the results of this paper.

\subsection*{Acknowledgments} The author sincerely thanks the anonymous referees for the helpful comments and valuable suggestions, which considerably improved the exposition of this work. 

\bibliographystyle{siamplain}

\bibliography{szhang}

\end{document}